\documentclass[10pt,article]{amsbook}
\usepackage[latin1]{inputenc}
\usepackage{amsmath,amssymb}

\newtheorem{theorem}{Theorem\ }
[section]\newcommand{\fdem}{\rightline{$\Box$}}
\newtheorem{lemma}[theorem]{Lemma\ }

\newtheorem{defi}[theorem]{Definition\ }
\newtheorem{fact}[theorem]{Fact\ }
\newtheorem{proposition}[theorem]{Proposition\ }

\newtheorem{notations}[theorem]{Notations\ }
\newtheorem{notation}[theorem]{Notation \ }

\newtheorem{property}[theorem]{Property\ }
\newtheorem{remark}[theorem]{Remark\ }

\def \C{\mathbb{C}}
\def \N{\mathbb{N}}

\def \H{\mathbb{H}}
\def \Z{\mathbb{Z}}
\def \R{\mathbb{R}}

\def \e{\epsilon}

\def \cH{\mathcal H}

\def \cP{\mathcal P}

\catcode`\?=\active\def?{\"u}\catcode`\?=\active\def?{\'e}
\catcode`\?=\active\def?{\^a}\catcode`\?=\active\def?{\`a}\catcode`\=\active\def{\c c}
\catcode`\=\active\def{\^e}\catcode`\=\active\def{\`e}\catcode`\?=
\active\def?{\"\i }\catcode`\?=\active\def?{\^\i }\catcode`\?=\active\def?{\^o}\catcode
`\?=\active\def?{\^u}\catcode`\=\active\def{\`u}\catcode`\¶=\active\def¶{\delta }
\catcode`\¹=\active\def¹{\pi }

\def\ds{\displaystyle}

\def \fdem{$\Box$}

\begin{document}
\centerline{\Large On some exotic Schottky groups}
\vspace{2mm}

 \centerline{Marc Peign\'e  
 $^($\footnote{Marc Peign\'e   LMPT, UMR 6083, Facult\'e des Sciences etTechniques, Parc de Grandmont, 37200 Tours -- mail : peigne@lmpt.univ-tours.fr}$^)$}
\vspace{2mm}

\centerline{May 2010}

\bigskip

 {\small {\bf Abstract. } We contruct  a Cartan-Hadamard manifold with pinched negative curvature whose group of isometries possesses    divergent discrete free subgroups  with parabolic elements   who do not satisfy the  so-called ``parabolic gap condition'' introduced in  \cite{DOP}. This construction relies on the comparaison  between the Poincar\'e series of these free  groups and the potential of some  transfer operator which appears naturally in this context. }

\section{Introduction}
        
 Throughout   this paper, $X$ will denote a  complete and simply connected Riemannian
manifold of dimension $ N\geq 2$  whose sectional curvature  is bounded between two negative constants  $-B^2\leq -A^2<0$. 
We denote by $d$ the distance on $X$ induced by the Riemannian metric and by $\partial X$ the boundary at infinity  ; the isometries of $X$ act as conformal transformations on $\partial X$  when it is endowed   by the so-called Gromov-Bourdon metric.

A  {\it Kleinian group} of $X$ is a non elementary  torsion free and discrete subgroup $\Gamma$ of orientation preserving isometries of $X$ ; this group  $\Gamma $ acts freely and properly discontinuously  on  $X$ and the quotient manifold  $M:=  X/\Gamma$ has a  fundamental group which can be identified with  $\Gamma$.   One says that $\Gamma$ is {\it a  lattice } when the Riemannian volume of $X/\Gamma$ is finite.

The limit set $\Lambda_\Gamma$ of a Kleinian group $\Gamma$ is the least non empty  $\Gamma$-invariant subset of $\partial X$ ; this is also the set of accumulation points of some (any) orbit $\Gamma\cdot{\bf x}$ of ${\bf x} \in X$. This set is of interest for further reasons ; in particular,  if $(\phi _t)_t$ denotes the geodesic flow on the unit tangent 
bundle  $T^1(X/\Gamma)$ of $X/\Gamma$, its {\it  non-wandering set}
$\Omega_\Gamma$  coincides with the projection  on   $T^1(X/\Gamma)$ of the set of unit 
tangent vectors on $X$ whose  points at infinity in both directions belong to
 $\Lambda_\Gamma$.

Note that the convex-hull $C(\Lambda_\Gamma)$ of $\Lambda_\Gamma$   is a $\Gamma$-invariant closed subset of $X$ and that the projection of $\Omega_\Gamma$ onto the manifold $X/\Gamma$ is in fact equal to $C(\Lambda_\Gamma)/\Gamma$.The group $\Gamma$ is said {\it convex 
cocompact} when it acts co-compactly on 
$C(\Lambda_\Gamma)$ and more generally {\it geometrically finite} when it acts like a lattice on some (any) $\epsilon$-neigbourhood  $C^\epsilon(\Lambda_\Gamma)$ of $C(\Lambda_\Gamma)$ (in otherwords, when  $vol( C^\epsilon(\Lambda_\Gamma)/\Gamma)<+\infty$ for some (any) $\epsilon >0$).

It is shown in  \cite {OP} that the existence and unicity of a measure of maximal entropy for the geodesic flow restricted to $\Omega_\Gamma$ 
    is  equivalent to the finiteness of a natural invariant Radon measure  on $T^1(X/\Gamma)$ with support $\Omega_\Gamma$,  
the so-called {\it Patterson-Sullivan measure} $m_\Gamma$. 
 In this paper, we  construct  
examples of isometry groups $\Gamma$ for which the restriction of the geodesic flow 
$(\phi _t)_t$ to the set $ \Omega_\Gamma $ exhibits  particular properties with 
respect to ergodic theory.   In particular,  for those  groups,  the Patterson-Sullivan measure   may  be  
 infinite and the associated dynamical
system $(\phi _t,\Omega_\Gamma )$ will thus have no measure  of maximal  entropy.

 We  now recall briefly the construction of the Patterson-Sullivan  measure associated with a Kleinian group $\Gamma$. The critical exponent of $\Gamma$ is the exponential growth of its orbital function defined by  
 $$
\displaystyle{\delta_\Gamma := \limsup _{r\to \infty} \frac{1}{r} 
\log 
\mbox{\rm  card}
\{\gamma \in \Gamma/ d({\bf x},\gamma \cdot {\bf x})\leq r \}}.
$$ 
It does not depend  on ${\bf x} \in X$ and coincides with the exponent of convergence of the    Poincar\'e series     of  $\Gamma$  defined by
${\bf P}_\Gamma(s, {\bf x}):= \sum _{\gamma \in \Gamma }e^{-sd({\bf 
o},\gamma \cdot {\bf o})}$ 
 ; this  series converges if $s>\delta_\Gamma$ and diverges when $s<\delta_\Gamma.$
The group $\Gamma$ 
is {\it divergent}   when the Poincar\'e series diverges at the 
critical exponent ; 
otherwise  $\Gamma$ is   {\it convergent}.

A construction due to Patterson  in constant curvature provides a family of $\delta_\Gamma$-{\it conformal measures}  
$\sigma=(\sigma_{\bf x})_{{\bf x}\in X }$
supported on the limit set 
$\Lambda_\Gamma$. D. Sullivan showed  also how to assign to $\sigma$  
an invariant measure 
for
the geodesic flow  $(\phi_t)_t$ restricted to $\Omega_\Gamma$. 
This construction has been extended by several people 
to the situation of a variable curvature space 
$X$ and an arbitrary  Kleinian  group $\Gamma$
acting  on it  \cite{K1}, \cite{Y}.

 It is important to recall that   the family of measures $\sigma$
associated with $\Gamma$ is  unique if and only if $\Gamma$ is divergent (see  \cite{Ro} for a complete statement). In this case, the corresponding $(\phi_t)_t$-invariant measure    constructed by Sullivan  depends only  on $\Gamma$, it is    the Patterson-Sullivan measure $m_\Gamma$ of $\Gamma$.

  We review now some basic results concerning the finiteness 
of  the measure $m_\Gamma$.
When $\Gamma$ is convex-coccompact, this measure is  of course finite since it is a Radon measure with compact support. The same property 
holds  when $\Gamma$ is a geometrically finite group acting on a
    locally symmetric space \cite {Su2}, \cite{CI} ; nevertheless, there exist non-geometrically finite groups with finite Bowen-Margulis measure \cite{P}.

The situation is much more complicated
in the general  variable curvature case, even for geometrically finite 
groups, because of the existence of parabolic subgroups. 

  There exist in particular   criteria which ensure  that a geometrically finite group $\Gamma$
 is divergent, for instance when
its Poincar\'e exponent $\delta_\Gamma$ is  strictly greater than the one of each of 
its parabolic subgroups \cite [Th\'eor\`eme A]{DOP}. This is the so-called {\it  
parabolic gap condition} (PGC), which is satisfied in particular 
when  the  parabolic subgroups of $\Gamma$ are themselves  
divergent. Furthermore,    the Patterson-Sullivan measure $m_\Gamma$
of a divergent geometrically finite group $\Gamma$ is finite if and 
only if, for 
any parabolic subgroup $ {\mathcal P}$ of $\Gamma$,  one has 
\begin{equation}\label{BMfinie}
\displaystyle{\sum_{p\in {\mathcal P}}d({\bf o}, p \cdot {\bf o})
e^{-\delta_{\mathcal P}d({\bf o}, p \cdot {\bf o})}<+\infty},
\end{equation}
 where  
 $\delta_{\mathcal  P}$ denotes the critical exponent of $\mathcal P$ 
\cite[Th\'eor\`eme B] {DOP} ;  this holds  in particular
  when the critical gap 
 property is satisfied and in this case, by the Poincar\'e recurrence theorem,    the geodesic flow $(\phi_t)_t$ is completely conservative with respect to $m_\Gamma$.

When $\Gamma$ is convergent, the Patterson-Sullivan
 measure $m_\Gamma$ 
 is infinite  and the geodesic flow $(\phi_t)_t$ is completely dissipative with respect to $m_\Gamma$.
 On may choose the metric in such a way there exist  
non elementary geometrically finite groups of 
convergent type ; in this case, the parabolic gap condition  is not 
satisfied and    the parabolic subgroups of $\Gamma$ of maximal 
Poincar\'e exponent  are convergent. 
In \cite{DOP} an explicit construction of such groups  is proposed.

As far as we know, there were no 
  examples of   geometrically finite 
groups $\Gamma$ of divergent type which do not satisfy 
the critical gap property  ; this contruction 
is really of 
interest because it gives   examples 
of geometrically finite manifolds for which 
the geodesic flow is completely 
conservative with respect to $m_\Gamma$ but
this measure  is infinite.   
We have the

\begin{theorem}\label{maintheorem}
  There exist Hadamard manifolds   with pinched negative 
 curvature  whose group of isometries contains geometrically finite Schottky groups $\Gamma$ of 
divergent  type which do not satisfy the
 parabolic gap condition PGC. 
 Furthermore, the    Patterson-Sullivan measure $m_\Gamma$ may be finite or infinite.
\end{theorem}

The paper is organized as follows : 
Section \S2 deals with  the construction of convergent parabolic groups  ; we recall in particular the results presented in \cite{DOP}.  
Section \S3 is devoted to the construction of Hadamard manifolds  containing convergent parabolic elements  and whose groupe of isometries is non elementary. 
In section \S4 we construct Schottky groups with convergent parabolic factor and  we explain how to choose  the metric inside the corresponding cuspidal end  to prove   Theorem \ref{maintheorem}.

 We fix here once and for all   some notation  about asymptotic behavior of functions :

\begin{notations} Let $f, g$ be two functions from $\R^+$ to $\R^+$.
We shall write
$f\stackrel{c}{\preceq}g $ (or simply $f \preceq g $) when $f(R) \leq c g( R)$ for some constant $c>0$ and $R$ large enough.
The notation
$f \stackrel{c}{\asymp}g $ (or simply $f \asymp g $) means  $f \stackrel{c}{\preceq}g \stackrel{c}{\preceq}f .$

 Analogously, we whall write $f \stackrel{c}{\sim}g $ (or simply $f \sim g $) when $|f(R) -  g( R)|\leq c$ for some constant $c>0$ and $R$ large enough.
\end{notations}

\section{On the existence of convergent parabolic  groups}

\subsection{The real hyperbolic space}
 We first consider the real hyperbolic space of dimension $N \geq 2$, identified  to the upper half-space ${\mathbb H}^N:=\R^{N-1}\times \R^{*+}$. In this model, the Riemannian hyperbolic metric   is given by 
 $\displaystyle {dx_{ }^2+dy^2\over y^2}$ where $dx_{ }^2+dy^2$ is the classical  euclidean metric on $\R^{N-1}\times \R^{*+}$. We denote by ${\bf i}$ the origin $(0, \cdots, 0, 1)$ of $\H^N$ and by $\Vert .\Vert$ the euclidean norm in $\R^N$.
 
 Let $p$ be a parabolic isometry of $\H^N$ fixing $\infty$ ; its induces on $\R^{N-1}$ an   euclidean isometry  which can be decomposed as  the product $p= R_p\circ T_p=T_p\circ R_p$ of an affine rotation $R_p$ and a translation $T_p$  with vector of translation $\vec{s}_p$. By an elementary calculous in hyperbolic geometry, one may check that the  sequence $(d({\bf i}, p^n\cdot {\bf i})-2\ln n\Vert \vec{s}_p\Vert)_{n \geq 1}$  converges to $0$.
 The Poincar\'e exponent of the group $\langle p\rangle$ is thus equal to ${1\over 2}$ and   $\langle p\rangle$ is divergent.

 More generally, for any parabolic subgroup ${\mathcal P}$ of the group of isometries of $\H^N$, the sequence $(d({\bf i}, p \cdot {\bf i})-2\ln  \Vert \vec{s}_p\Vert_{ })_p$
 converges to 0 as $p\to  \infty$ in ${\mathcal P}$. By one of Bieberbach's theorems, the group ${\mathcal P}$  contains a finite index  abelian subgroup  ${\mathcal Q}$  which acts by translations on a subspace $\R^k$ of $\R^{N-1}$ ; in other words, there exist $k$ linearly independant vectors $\vec{s}_1, \cdots,
 \vec{s}_k$ and a finite  set $F\subset {\mathcal P}$ such that  any $p \in \mathcal P$ may be decomposed as $p=p_{\vec{s}_1}^{n_1} \cdots p_{\vec{s}_k}^{n_k} f$ with $n_1, \cdots , n_k \in \Z$ and $f\in F$ so that 
\begin{eqnarray*}
 {\bf P}_{\mathcal P}(s) =1+\sum_{p \in \mathcal P^*}  e^{-sd({\bf i}, p\cdot{\bf i})}&=&1+ \sum_{p \in \mathcal P^*} {e^{s o(p)}\over \Vert \vec{s}_p\Vert^{2s}}\\
 &=&1+ \sum_{f \in F}\ \ \sum_{\bar n=(n_1, \cdots , n_k)\in (\Z^k)^*}
{e^{so(\overline{n})}\over 
\Vert n_1\vec{s}_1+ \cdots +n_k\vec{s}_k\Vert^{2s}
}. \end{eqnarray*}
The Poincar\'e exponent of $\mathcal P$ is thus equal to ${k\over 2}$ and   the group  is divergent.

 All these calculous may be done in the following (less classical) model : using the natural diffeomorphism between $\H^N$ and $\R^N$ defined by $(x, y)  \mapsto (x, t):= (x, \ln y)$
one may  endow $\R^N$ with the hyperbolic metric $\displaystyle g_{hyp}:= e^{-2t} dx_{ }^2+dt^2$. 

In this model, we fix the origin  ${\bf o}=(0, \cdots, 0)$ and the vertical lines $\{(x, t)/t \in \R\}$ are clearly geodesics. For any $t \in \R$, we denote by ${\mathcal H}_t$ the hyperplane $\{(x, t): x\in \R^{N-1}\}$  ; this corresponds  to  the horosphere centered at $+\infty$ and passing through $(0,\cdots, 0, t)$. For any $x, y \in \R^{N-1}$, the  distance between ${\bf x}_t:=(x, t)$
 and ${\bf y}_t:= (y, t)$ for the   metric $ e^{-2t} dx_{ }^2$ induced by $g_{hyp}$ on ${\mathcal H}_t$ is equal to $e^{-t}\Vert x-y\Vert_{ }$ ; furthermore, if $t$ is choosen in such a way that this distance is equal to $1$ (namely $t=\ln \Vert x-y\Vert_{ })$, then the union of the 3 segments $[{\bf x}_0, {\bf x}_t],
 [{\bf x}_t, {\bf y}_t] $ and $[{\bf y}_t, {\bf y}_0]$ lies at a bounded distance of the hyperbolic geodesic joigning  ${\bf x}_0$ and ${\bf y}_0$ which readily implies that $d({\bf x}_0, {\bf y}_0)-2\ln \Vert x-y\Vert_{ }$ is bounded. 
 
 This crucial fact is the key to understand geometrically the estimations above ; it first appeared  in  \cite{DOP} and allowed the authors to construct negatively curved manifolds with convergent parabolic subgroups, we recall in the following subsection this construction.
 
 \subsection{ The metrics $T_{ a, u}$ on $\R^N$}

 We consider on $\R^{N-1}\times \R$   a
Riemannian metric   of the form
$g=T^2(t)dx_{ }^2+dt ^2$,
where $dx_{ }^2$ is a fixed   euclidean  metric 
on  $\R^{N-1}$ and  $T: \R\to \R^{*+}$ is
a $C^{\infty}$   non-increasing function.
The group of isometries of  $g$ contains the isometries of 
$\R^{N-1}\times \R$   fixing the last  coordinate. The sectionnal 
curvature 
at  $\displaystyle (x,t)=(x_1, ...,x_{N-1},t)$ does not depend on $x$ 
: it is  $\displaystyle K(t)=-\frac 
{T''(t)}{T(t)}$ on any   plane  
$\displaystyle \langle \frac {\partial }{\partial X_i}, \frac 
{\partial}{\partial t}\rangle, 1\leq i\leq N-1$, 
and  $-K^2(t)$ on any plane  
$\displaystyle \langle \frac {\partial }{\partial X_i}, \frac 
{\partial}{\partial X_j}\rangle, 1\leq i<j\leq N-1$ (when $N \geq 2)$.

   It is convenient to consider the 
  non-decreasing function
$u : \R^{*+}\to \R$ satisfying the following implicit equation 
 \begin{equation}\label{implicit}
 T(u(s))=\frac 1s.
 \end{equation}
Then, the value of the  curvature of $g$ is :
\begin{equation}\label{curvature}
K(u(s)):=-\frac {T''(u(s))}{T(u(s))}=
-\frac {2 u'(s)+s u''(s)}{s^2(u'(s))^3}.
\end{equation} 
Note that  $g$ has negative curvature if and only if $T$ is 
convex.    
For instance, we have seen in the previous subsection that for $\displaystyle u(s)=  \log s$ one gets $T(t)= e^{-t}$ and obtains a 
model of the hyperbolic space   of constant curvature $-1$.

      As it was seen in \cite{DOP}, the function $u$ is of interest since it gives precise estimates (up a bounded term) of the distance between points lying on the same  horosphere  ${\mathcal H}_t:= \{(x, t): x \in \R^{N-1}\}$ where $t \in \R$ is fixed.  
Namely, the  distance between ${\bf x}_t:=(x, t)$
 and ${\bf y}_t:= (y, t)$ for the   metric $ T^2(t)dx_{ }^2$ induced by $g$ on ${\mathcal H}_t$ is equal to $T(t)\Vert x-y\Vert_{ }$ ; for $t=u( \Vert x-y\Vert_{ })$, this distance is thus  equal to $1$, and the union of the 3 segments $[{\bf x}_0, {\bf x}_t],
 [{\bf x}_t, {\bf y}_t] $ and $[{\bf y}_t, {\bf y}_0]$ lies at a bounded distance of the hyperbolic geodesic joigning  ${\bf x}_0$ and ${\bf y}_0$ (see  \cite{DOP}, lemme 4)  : this readily implies that $d({\bf x}_0, {\bf y}_0)-2u( \Vert x-y\Vert)$ is bounded.
 
 In the sequel, we will assume that the function $u$ coincides with the function $s \mapsto \ln s$ on $]0, 1]$ ; in otherwords, the  restriction to the set $]0, 1]$ of the corresponding function $T_u(t)$ satisfying (\ref{implicit}) is equal to $t \mapsto e^{-t}$. More generally, we will ``enlarge '' the  area where $T_u(t)$ and $e^{-t}$ coincides to the domain $\R^{N-1}\times ]-\infty, a]$ with $a$ arbitrary,  introducing the following 
 \begin{notation}\label{metriquedanscusps} Let $a\in \R$ and $u: \R^{*+} \to \R$ be a $C^2$ non decreasing function  such that 
 \begin{itemize}
 \item $\ u(s)=\ln s$  for any $s \in ]0, 1]$
 \item
 $K(u(s)) \in [-B^2, -A^2]\subset \R^{*-}$ for any $s>0$.\end{itemize}
We  endow  $\R^{N-1}\times \R$  with the metric  
$
T_{a, u}^2(t) dx^2+dt^2, 
$
where $T_{a, u}$ is given  by  
\begin{equation}\label{T_{a,u}}
\forall t \in \R \quad T_{a,u}(t) := \Bigl\{
\begin{matrix}e^{-t} & if &t\leq a \\ 
\displaystyle {e^{-a}\over u^{-1}(t-a)} & if &t\geq a \
\end{matrix}
\end{equation}.
\end{notation}
 Note that this metric    has constant curvature  $-1$ on the domain $\R^{N-1}\times ]-\infty, a].$

\subsection{On the existence of metrics with convergent parabolic  groups} 
In this paragraph, we fix $a\in \R$ and endow $\R^{N-1}\times \R$ with the  metric $T_{a,  u}^2(t) dx^2+dt^2$ where  $\displaystyle u(s)=     {\ln s+ \alpha  \ln \ln  s}  $ for $s$ large enough  and some constant $\alpha >0$ ; in this case,  the curvature   varies, nevertheless  one 
has 
$\displaystyle \lim _{s\rightarrow  \infty}K(u(s))=-1$ and
all derivatives of 
$K(u(s))$ tend    to $0$ as $s \to +\infty$ . We  will  first need the following
\begin{lemma} \label{u_alphaexiste} Fix $\kappa \in ]0, 1[.$ 
For any $\alpha \geq 0$, there exists  a constant $s_\alpha \geq 1$  and a  non decreasing $C^2$ function $u_\alpha: \R^{*+}\to \R$ such that 
  \begin{itemize}
 \item $u_{ \alpha}(s)=\ln s$ \mbox{\rm if}
$0<s\leq 1$ 
\item $u_{ \alpha}(s)= \ln s+ \alpha \ln  \ln  s $  \mbox{\rm if} 
$s \geq s_\alpha$.
\item $\displaystyle K(u_\alpha(s)):=
-\frac {2 u_\alpha'(s)+s u_\alpha''(s)}{s^2(u_\alpha'(s))^3} \leq -\kappa^2.$
\end{itemize}
\end{lemma}
 \noindent Proof.  We first fix a $C^2$ non decreasing function $\phi: \R\to [0, \alpha]$, which vanishes on $\R^{-}$ and is equal to $\alpha$ on $[1, +\infty[$. For any $\epsilon >0$, we consider the function $v_\epsilon : [e, +\infty[\to \R$ defined by
$$
\forall s \geq 1 \quad v_\epsilon(s):= \ln s+\phi_\epsilon(s)\ln \ln s
$$
where $\phi_\epsilon(s):= \phi(  \epsilon  \ln \ln s)$.
 A straightforward computation gives, for any $s \geq e$ 
 $$
 \frac {2 v_\epsilon'(s)+s v_\epsilon''(s)}{s^2(v_\epsilon'(s))^3}= {N_\epsilon(s)\over D_\epsilon(s)}
 $$
 with \begin{itemize}
 \item 
 $\displaystyle N_\epsilon(s):= 1+{\phi_\epsilon(s)\over \ln s}-{\phi_\epsilon(s)\over (\ln s)^2}
 +2 \phi'_\epsilon(s) \Bigl(s\ln \ln s +{s \over \ln s}\Bigr)+ \phi_\epsilon''(s)s^2\ln \ln s$, 
 \item
 $\displaystyle D_\epsilon(s):= \Bigl(1+{\phi_\epsilon(s)\over \ln s}
 + \phi'_\epsilon(s) s\ln \ln s\Bigr)^3$,
 \item $\displaystyle\phi_\epsilon'(s)={\epsilon \over  s \ln s}
 \phi_\epsilon(s) $ and $\displaystyle\phi_\epsilon''(s)={\epsilon^2-\epsilon(1+\ln s)\over(s\ln s)^2} \phi_\epsilon(s).$
 \end{itemize}
 For any continuous function $g : [e, +\infty[$  converging to $0$ at infinity, one gets  $g \phi_\epsilon \to 0$ uniformly on $[e, +\infty[$ ; consequently, one obtains,   as $\epsilon\to 0$ and uniformly  on $[e, +\infty[$ 
 $$
 \phi'_\epsilon(s) \Bigl(s\ln \ln s +{s \over \ln s}\Bigr)=\epsilon\Bigl(
 {\ln \ln s\over \ln s}+{1\over (\ln s)^2}
 \Bigr)\phi_\epsilon(s) \to 0
 $$
 and 
$$
 \phi_\epsilon''(s)s^2\ln \ln s={\ln \ln  s\over (\ln s)^2}\Bigl(\epsilon^2-\epsilon(1+\ln s)\Bigr)\phi_\epsilon(s) \to 0,
$$ 
  so that 
 $
\displaystyle  \frac {2 v_\epsilon'(s)+s v_\epsilon''(s)}{s^2(v_\epsilon'(s))^3}\to 1$. One may thus choose $\epsilon_0 >0$ such that 
 $$
\forall s \geq e \qquad - \frac {2 v_{\epsilon_0}'(s)+s v_{\epsilon_0}''(s)}{s^2(v_{\epsilon_0}'(s))^3}\leq -\kappa^2
 $$
 and one sets
 \begin{eqnarray*}
u_\alpha(s) := \Bigl\{
\begin{matrix}
\ln s & {\rm if}  & 0<s\leq e\\ 
v_{\epsilon_0}(s) & {\rm if} & s\geq e,
\end{matrix}
\end{eqnarray*}
with $s_\alpha:=\exp(\exp(1/\epsilon_0))$.
\fdem

We  thus fix $a, \alpha \geq 0$ and 
endow $\R^N=\R^{N-1}\times \R$ with the metric $ T ^2_{a, u_\alpha}(t)dx^2+dt^2 $ 
where $u_\alpha$ is given by   Lemma \ref{u_alphaexiste}.
This metric   has pinched negative curvature less than $-\kappa^2$ and constant negative curvature in the domain $\{(x, t): t\leq\  a\}$.

  Now, let $\cP $ be a discrete  group of isometries  of
$  \R^{N-1}$
of
rank $k \in \{1, \cdots, N-1\}$, i.e generated by $k$ linearly independent translations $p_{\vec{\tau}_1}, \cdots, p_{\vec{\tau}_k}$ in $  \R^{N-1}$.   In order to simplify the notations, $\bar n=(n_1, \cdots, n_k) \in \Z^k$ will represent the translation of vector $ n_1\vec{\tau}_1+ \cdots +n_k\vec{\tau}_k$ and  $\vert \bar n \vert$ will denote its euclidean norm. These translations are also isometries of $\R^N$ endowed with the  metric $T_{a, u_\alpha}(t)^2dx^2+dt^2$ given above and the  corresponding Poincar\'e series of $\cP$  is given by
\begin{eqnarray*}
 {\bf P}_{\mathcal P}(s) =1+\sum_{p \in \mathcal P^*}  e^{-sd({\bf o}, p\cdot{\bf o})}&=&1+ \sum_{\bar n \in (\Z^k)^*}{e^{-2su_{a, \alpha}(\vert \bar n\vert)-s O(\bar n)}}\\
 &=&1 + \sum_{\bar n \in (\Z^k)^*}
{e^{-sO(\overline{n})}\over 
\vert \bar n \vert^{2s}
\Bigl(\ln \vert \bar n \vert\Bigr)^{2s\alpha}
}. \end{eqnarray*}
We have thus prove the 
\begin{proposition}
Let $\R^N$ be endowed with  the metric $ T ^2_{a, u_\alpha}(t)dx^2+dt^2 $ 
where $u_\alpha$ is given by   Lemma \ref{u_alphaexiste}.  If $\cP$ is a discrete group of isometries of $\R^{N-1}$ of rank $k$, its critical  Poincar\'e exponent   is equal to ${k/2}$ ; furthermore, the group $\cP$ is convergent if and only if $\alpha >1$.
\end{proposition}  

\begin{remark} One may also choose $u$ is such  a way that 
$u^{-1}(t) = e^{t/2-\sqrt{t}}$. If $r=1$, the critical exponent of 
the 
associated Poincar\'e series is equal  to ${1\over 2}$ and the group $\cP $ 
is 
also convergent ; this last example appears in \cite{Sch},  where    
some 
explicit results are 
given,  in terms of the Poincar\'e series of the parabolic 
groups,  which guarantee the equidistribution of  the horocycles  on 
geometrically finite  negatively curved surfaces.
\end{remark}

\section{Weakly homogeneous Hadamard manifolds of type $(a, u_\alpha )$}

In the previous section,  we have endowed $\R^N$ with a metric   $T_{a, u }(t)^2dx^2+dt^2$  ; unfortunately,  in this construction, excepted for some particular choice of $u$,  all  the isometries   fix the same point at infinity and the group  $Is(\R^N)$ is thus elementary.
We need now to construct  an Hadamard manifold with a metric of  this inhomogeneous   type  in the neighbourhood of some points at infinity but whose group of isometries is non elementary.

\subsection{Metric  of type $(a, u )$ relatively to some group $\Gamma$ and some horoball $\mathcal H$}
Consider first  a  non uniform lattice   $\Gamma$ of isometries of 
$\H^N$. The manifold $M:=\H^N/\Gamma$  has finite volume but is not 
compact ; it thus possesses finitely many cusp 
 $C_1, \cdots, C_l$, each cusp $C_i$ being  isometric to the quotient  of 
some horoball  ${\mathcal H}_i$ of $\H^N$ (centered at a point $\xi_i$) by a  Bieberbach group
$\cP_i$ with rank $N-1$. Each group $\cP_i$  also acts by isometries on $\R^{N-1}\times \R$ endowed 
 with one of the metrics  $ T_{a, u}(t)^2 dx^2+dt^2$ given by Notation \ref{metriquedanscusps}. 
 
Now, we endow  $\R^{N-1}\times \R$ with one of these metrics $T_{a, u}(t)^ 2dx^2+dt^2$   and  choose $a$ in such a way  we may paste  the quotient
 $(\R^{N-1} \times [0, +\infty[)/\cP_1 $ with 
$M \setminus C_1$. The Riemannian manifold $M$ remains   negatively curved with finite volume. By construction,  the group $\Gamma$ acts isometrically on the
 universal covering $X\simeq \R^N$ of $M$ endowed with the lifted metric $  g_{a, u}$  ; note that $g_{a, u}$ coincides  with the metric $T^2_{a, u }(t)dx^2+dt^2$  on the preimage by $\Gamma$ of the cuspidal end $C_1  \ ^($\footnote{By the choice of $C_1$, the horoballs $\gamma\cdot \cH_1, \gamma \in \Gamma,$ are disjoint or coincide, they are also isometric to $\R^{N-1}\times \R^+$ endowed with the hyperbolic metric $e^{-2t}dx^2+dt^2$ ; another way to endow $\R^N$ with the new metric  $g_{a, u }$ is to replace inside each horoball $\gamma\cdot \cH_1$ the hyperbolic metric with the restriction of $T^2_{a, u }(t)dx^2+dt^2$   to the half space $\R^{N-1}\times \R^+$}$^)$.

All this discussion gives sense to the following definition    :
 \begin{defi}\label{Hadamardtypealphaa} Fix $ a, \alpha \geq 0$, let $u_\alpha$ be the function given by Lemma \ref{u_alphaexiste} and  $(X, g)$  a negatively curved Hadamard manifold whose group of isometries contains a non uniform lattice $\Gamma$.
 
 Assume that $X/\Gamma$ has one cusp $C$, let $\cP$ be  a maximal parabolic subgroup  of $\Gamma$ corresponding to this cusp, with  fixed point  $\xi $ and let $\cH $ be an   horoball   centered  at $\xi$ such that   the $\gamma\cdot \mathcal H, \gamma \in \Gamma$, are disjoints or coincide.

One endows the manifold $X$  with the metric $g_{a, u}$ defined by
\begin{enumerate}
\item  $g_{a, u}$ has constant curvature $-1$ outside the set  
$ \bigcup_{\gamma \in \Gamma}\gamma\cdot \cH$ 
\item $g_{a, u}$ coincides 
with the metric  
$T_{a, u}(t)^2dx^2+dt^2$  inside each horoball $\gamma\cdot \cH, \gamma \in \Gamma$.
\end{enumerate}

 One says that {\bf  the Riemannian manifold   $(X, g_{a, u})$  has type  $( a, u)$     relatively to  the group $\Gamma$ and the horoball $\cH$.}
More generally, one says that  {\bf $(X, g)$ has type $u$ } when, for some $a \in \R$,  some lattice    $\Gamma$ and some horoball $\cH$,  it has type  $( a, u)$  relatively to    $\Gamma$ and  $\cH$.

   \end{defi}
   
   \begin{remark}\label{stripe}
If  the metric $g$  has    type  $( a, u )$   relatively to  $\Gamma$ and  $\cH$,   the curvature remains equal to $-1$ in the stripe $\R^{N-1}\times [0, a]\subset  \cH$. In the limit case ``$a=+\infty$'', one refinds the hyperbolic metric of constant curvature $-1$. 
  \end{remark}

  By construction, the elements of 
$\Gamma$ are isometries of $(X, g_{a, u})$.
 It is a classical fact that the group of isometries of $\H^N$ is quite large since in particular it acts transitively on the hyperbolic space (and even on its unit tangent bundle). This property remains valid when $X$ is symmetric, otherwise its isometry group is discrete (\cite{E}, Corollary 9.2.2) ; consequently, if $g$ has type $( a, u_\alpha)$ with $\alpha >0$, the group of isometries of  $(X, g)$    does not inherit this property  of transitivity,  it is  is discrete and  has $\Gamma$ as   finite index  subgroup (\cite{E}, Corollary 1.9.34).

  We fix now one and for once the following 
  
   \begin{notation}\label{NOTATION} From now on, we  consider an Hadamard manifold  $X$, with  origin ${\bf o}$, whose group of isometries contains a non uniform lattice $\Gamma$ ;  we fix  a maximal parabolic subgroup  $\cP$ of $\Gamma$  with fixed point  $\xi \in \partial X$ and  an  horoball   $\cH$ centered at $\xi $ such that the horoballs $\gamma\cdot \mathcal H, \gamma \in \Gamma$, are disjoints or coincide.
 
 {\bf We fix $\alpha \geq 0$ and we assume that, for any $a\geq 0$, the manifold  $X$ may be endowed 
 with a metric $g_a:= g_{a,u_\alpha}$ of type $(a, u_\alpha)$  relatively to $\Gamma$ and $\cH$}, where $u_\alpha$ is given by Lemma \ref{u_alphaexiste}.

We denote   by $d_a$ the corresponding distance on $X$.
\end{notation}
 Note that, by construction, the sectional curvature of  
$g_{a}$ is 
pinched between two non positive constants and is less than $ -\kappa^2$  for somme constant $\kappa > 0$ which does not depend on 
$a$. Furthermore, using the  fact that  $u_\alpha$ is non negative and   uniformly continuous on $[1, +\infty[$, one obtains the
  \begin{property}\label{quasi-isometry} 
For any $a, a'$ and $\alpha \geq 0$, there exists a constant $K = K_{a, a', \alpha} \geq 1$ with $K_{a, a', \alpha}
 \to 1$ as $a'\to a$ such that
  $$\frac{1}{K}\ g_{a }\leq g_{a' }\leq  K  g_{a},$$
  so that
$$
 \frac{1}{K }\ d_{a}\leq d_{a'}\leq \ K d_{a}.
 $$
   \end{property}
\begin{remark}Note that if $a'>a$,  one has  in fact  $g_a\geq g_{a'}$ and so $d_a\geq d_{a'}$. It will be used in the last section.
\end{remark}

 \subsection{On the metric structure of the boundary at infinity}

In this paragrah we describe the metric structure of the boundary at infinity of   $X$; we  need first  to consider  the Busemann function  $ {\mathcal  
B}_{\cdot}^{(a)}(\cdot, \cdot)$ defined
by :  

for any $x \in \partial X$ and any ${\bf p}, {\bf q}$ in $X$
$${\mathcal 
 B}_{x}^{(a)}({\bf p},{\bf q}) = \lim_{{\bf x} \to x}d_{a}({\bf p}, 
{\bf x})
 -d_{a}({\bf q}, {\bf x}).
 $$ 
The Gromov product on $\partial X$, based at the origin ${\bf o}$,  between the points $x$ and $y$ in $\partial X$ is   defined by 
$$
(x\vert y)^{(a)}:= \frac{ {\mathcal  
B}_{x}^{(a)}({\bf o}, {\bf z})+
  {\mathcal B}_{y}^{(a)}({\bf o}, {\bf z}) }{2}
$$
where ${\bf z}$  is any point  on the geodesic $(x,y)$ (note that the value of $(x\vert y)^{(a)}$  does not depend on ${\bf z})$.
By \cite{Bou}, the function 
\begin{eqnarray*}
D_a: \partial X\times \partial X &\to& \R^+\\
(x, y)&\mapsto& D_a(x, y) := \Bigl\{
\begin{matrix}
\exp \Bigl( -\kappa(x\vert y)^{(a)}\Bigr)& if &x\neq y\\ 
0&if& x=y
\end{matrix}.
\end{eqnarray*}
is a distance on $\partial X$ ; furthermore, the cocycle property satisfied by the  Busemann functions readily implies that for any $x, y \in \partial X$ and  $\gamma \in \Gamma$
\begin{equation}\label{TAF1}
D_{a}(\gamma\cdot x, \gamma\cdot y)=
 \exp(-{\kappa \over 2} {\mathcal  B}_{x}^{(a)}({\bf o}, 
\gamma^{-1}.{\bf 
  o})) \exp(-{\kappa \over 2} {\mathcal  B}_{y}^{(a)}({\bf o}, 
\gamma^{-1}.{\bf 
  o}))D_{a}(x, y).
\end{equation}
In other words, $\gamma$ acts on $(\partial X, D_{a} )$ as a conformal transformation with coefficient of conformality
$$
|\gamma'(x)|_{a}= \exp(-\kappa {\mathcal  B}_{x}^{(a)}({\bf o}, 
\gamma^{-1}.{\bf 
  o}))$$ at the point $x$, since equality (\ref{TAF1}) may be rewrite 
  \begin{equation}\label{TAF2}
D_{a}(\gamma\cdot x, \gamma\cdot y)=
\sqrt{\vert \gamma'(x)\vert_{a}\vert 
\gamma'(y)\vert_{a}
  }D_{a}(x, y).
\end{equation}

 We will need to control the regularity  with respect to $a$ of the Busemann function 
 $x\mapsto {\mathcal  B}_{x}^{(a)}({\bf o}, z)$. By Property \ref{quasi-isometry}, the   
  spaces 
  $(X, d_{0})$ and 
  $(X, d_{a})$  are quasi-isometric and, for any $a_0>0$,  there exist a 
  constant $K_{0}\geq 1$ such that 
 \begin{equation} \label{K_0 quasi isometrie}
 \forall a \in [0, a_0]\quad \frac{1}{K_{0}}d_{0}\leq d_{a}\leq K_{0} d_{0}.
\end{equation}
 Note that, by \cite{GH},  one also gets 
  \begin{equation}\label{gromovquasi-iso}
\frac{1}{K_{a, a', \alpha}}(y|z)_{a}\leq (y|z)_{a'}\leq K_{a, a', \alpha}(y|z)_{a}.
 \end{equation}

  The corresponding distances $D_{a}$ on $\partial X$ are thus 
  H\"{o}lder equivalent ; more precisely, we have the 
  \begin{property} \label{holder-equivalent}
  For any $a_0>0$, there exists a real $\omega_0\in ]0, 1]$ such that, for all $a \in [0, a_0]$,  one gets
  $$
  D_{0}^{1/\omega_0}\leq D_{a}\leq 
 D_{0}^{\omega_0}.
  $$
  \end{property}

 The regularity of the Busemann 
function 
 $x \mapsto {\mathcal  B}_{x}^{(a)}({\bf o},{\bf  p})$ where ${\bf  
p}$ is a fixed point in 
 $X$ is given by the following Fact,   which precises  a result due to M. Bourdon.

\begin{fact}\label{Fact}  \cite{BP} 
  Let $E \subset \partial X$ and $F \subset X$ 
  two  sets  whose closure $\overline{E}$ and $\overline{F}$ 
  in $X\cup \partial X$ are disjoint.  Then the family of functions 
  $x \mapsto {\mathcal  B}_{x}^{(a)}({\bf o}, {\bf  p})$, 
  with ${\bf  p} \in F$,  is equi-Lipschitz continuous  on $E$ with 
respect 
to $D_{a}$. 

In particular,  for $a_0>0$ fixed, there exist  $\omega \in ]0, 1[$ and $C>0$ such that, for all  $a \in [0, a_0]$, one gets
\begin{equation}\label{equi-holder}
\forall x, y \in E, \forall {\bf p} \in F 
\quad 
\Bigl\vert
{\mathcal  B}_{x}^{(a)}({\bf o}, {\bf  p})-{\mathcal  B}_{y}^{(a)}({\bf o}, {\bf  p})
\Bigr\vert
\leq D_0(x, y)^{\omega}.
\end{equation}
\end{fact}

 \section {Divergent Schottky groups without  PGC}

 \subsection{On the existence of convergent Schottky groups when $\alpha >1$ }
 
 The fact that $\alpha >1$ ensures that any subgroup of $\cP$ is convergent. In \cite{DOP}, it is proved that $\Gamma$ possesses also non elementary subgroups
 of convergent type ; we first  recall this construction and precise the statement.  
 \begin{proposition} \label{Schottkygroupsconvergentdivergent} There exist Schottky subgroups $G$  of   $\Gamma$  and $a_0\geq 0$ such that 
 \begin{itemize}
\item $G$ has Poincar\'e exponent  ${1\over 2}$ and is convergent  on $(X, g_0)$ 
\item $G$ has Poincar\'e exponent  $>{1\over 2}$ and is divergent  on $(X, g_a)$  for $a\geq a_0$. 
\end{itemize}\end{proposition}
 
 Note that the group $G$ necessarily contains a parabolic element, otherwise it would be convex co-compact and thus of divergent type, whatever metric $g_a$ endows $X$.
 
  \noindent Proof. We first work in constant negative curvature $-1$ and fix a parabolic isometry  $p\in \cP$. Since $\Gamma$ is non elementary, there exists
  an hyperbolic isometry $q \in \Gamma$  whose fixed points are distinct from the one of  $p$. If necessary, one may shrink the horoball $\cH$  in such a way that  the projection of the axis of $h$ on the manifold $M=\H^N/\Gamma$  remains outside the cuspidal end $C\simeq \cH/\cP$ ; in others words,  one may fix  ${\bf o}$ on the axis of $q$ and assume that  for any $n \in \Z^* $ the geodesic segments $[{\bf o}, q^n\cdot{\bf o}]$ lie  outside the set $\displaystyle \bigcup_{\gamma\in \Gamma} \gamma\cdot \cH$ (so that in the area where the curvature is constant when $X$ will be endowed with the metric $g_a$).

  By the dynamic of the elements of $\H^N$ there exist 
   two compact sets ${\mathcal U}_{p}$ and  ${\mathcal 
U}_{q}$ in 
$X\cup \partial X $ as follows  :
\begin{enumerate}
\item ${\mathcal U}_{p}$ is a neigbourhood of the fixed point  $\xi_p$ of $p$ ;
\item ${\mathcal U}_{q}$ is a neigbourhood of the fixed points 
$\xi_q^+$ and $\xi_q^-$ of $q$ ;
\item there exists $\theta>0$ such that for any ${\bf  x}\in{\mathcal U}_{p}$ and ${\bf y}\in {\mathcal U}_{q}$  the angle $\widehat{{\bf x}\ {\bf o}\ {\bf y}}$ is greater than $\theta$ ;
\item for all $k \in \Z^*$ one has 
$$q^k\Bigl((X\cup \partial X)\setminus 
{\mathcal U}_{q}\Bigr)\subset {\mathcal U}_{q} \quad \mbox{\rm and}\quad 
p^k\Bigl((X\cup \partial X)\setminus 
{\mathcal U}_{p}\Bigr)\subset {\mathcal U}_{p}.$$ 
\end{enumerate}
By the Klein's tennis table lemma,   the group $\langle p, q\rangle$  generated by  $p$ and $q$ is 
free. Therefore  each element    $\gamma \in\langle p, q\rangle$, $\gamma\neq Id,$ may be   decomposed  in a unique way
  as  
a product $\alpha_{1}^{n_{1}}\alpha_2^{n_2}\ldots \alpha_{k}^{n_{k}}$ with $\alpha_{i}
\in \{p, q\}, n_{i}\in \Z^*$ and  $\alpha_{i}\neq \alpha_{i+1}$ ; the integer $k$ is the {\bf length} of $\gamma$ and $\alpha_k$ is its {\bf last letter}.

Let us now endow $X$ with the metric $g_0=g_{0, u_\alpha}$. If ${\bf x} \in {\mathcal U}_{p}$ and ${\bf y} \in {\mathcal U}_{q}$, the path 
which is the disjoint union of the geodesic ray $({\bf x},{\bf o}]$ and $[{\bf 
o},{\bf y})$ 
is a quasi-geodesic  in $(X, g_0)$; therefore  there  exists a constant $C>0$   
which only 
depends on 
the sets ${\mathcal U}_{p}$ and ${\mathcal U}_{q}$ and on the bounds 
on the 
curvature - that is, on the choice of the function $u_\alpha$ - such that
$d_{0}({\bf x}, {\bf y})\geq d_{0}({\bf x}, {\bf o})+d_{0}({\bf o}, {\bf y})-C$
    The Poincar\'e series of  this group  equals 
    $$
{\bf P}_{\langle p, h\rangle}(s) = 1+ \sum_{l\geq 1}\sum_{m_{i}, n_{i} \in \Z^*}e^{-sd_{\bf o}({\bf o}, p^{m_{1}}q^{n_{1}}\cdots p^{m_{l}}q^{n_{l}}\cdot {\bf o})}
.$$
 It 
follows
that 
$$
{\bf P}_{\langle p, q\rangle}(s) \leq 1+  \sum_{l\geq 1}\bigl(
e^{2sC} \ \sum_{m\in \Z^*}e^{-sd_{0}({\bf o}, p^m\cdot {\bf o})}\ 
\sum_{n\in \Z^*}e^{-sd_{0}({\bf o}, q^n\cdot {\bf o})}\ \bigr)^l.
$$
Recall that 
$\displaystyle{d_{0}({\bf o}, p^m\cdot {\bf o})= 
2 \ln m+2\alpha \ln \ln \vert m\vert  +\  a\  bounded\ term}$ ; since  $\alpha >1$, the series 
$\sum_{m\in \Z^*}e^{-sd_{0}({\bf o}, p^m\cdot {\bf o})}$ converges at its critical 
exponent 
$\delta_{\langle p\rangle}={1\over 2}$. We may now 
replace  $q$ by a sufficient
 large power $q^k$ in order to get 
$$
e^{C} \
 \sum_{m\in \Z^*}e^{-{1\over 2}  d_{0}({\bf o}, p^m\cdot {\bf o})}
  \sum_{n\in \Z^*}e^{-{1\over 2} d_{0}({\bf o}, q^{kn}\cdot {\bf 
o})}\ 
\ 
<1.$$
 It comes out that the critical exponent of 
the group  $G$ generated by  $p$  and $h:=q^k$
is less than ${1\over 2}$ and that ${\bf P}_G({1\over 2})<+\infty$ ; since 
 $p \in G$, one also gets
$\delta_G\geq \delta_{\langle p\rangle}={1\over 2}$. Finally $\delta_{G}={1\over 2}$ and 
$G$ is convergent, with respect to the metric $d_{0}$ on $X$.

Let us now prove that for $a$ large enough the group $G$ is divergent on $(X, d_a)$. 
By the
 triangular inequality one first gets
 \begin{eqnarray*} 
 \sum_{g\in G}e^{-s d_a({\bf o}, g\cdot {\bf o})}&\geq&  \sum_{l\geq 1}\sum_{n_{i}, m_{i} \in \Z^*}e^{-sd_{a}(o, p^{m_1}h^{n_{1}}\cdots p^{m_{l}}h^{n_{l}} \cdot {\bf o})}
\\
&\geq&  \sum_{l\geq 1}\Bigl(
\sum_{n\in \Z^*}e^{-sd_{a}({\bf o}, p^{m}\cdot {\bf o})}
\sum_{m\in \Z^*}e^{-sd_{a}({\bf o}, h^n\cdot {\bf o})}\Bigr)^l.
\end{eqnarray*}
Recall first that, when the curvature is constant and equal to $-1$ (that is 
to say 
"$a= +\infty$" in the definition of $g_a$), 
the quantity $d_{\H^N}({\bf o}, 
p^m\cdot {\bf o})-2\log \vert m\vert $ is bounded,  so  the parabolic group
$\langle p\rangle$ is
divergent with critical exponent   ${1\over 2}$.
There thus exists $\epsilon_0>0$ such that, for $\epsilon \in ]0, \epsilon_0]$, one gets 
\begin{equation}\label{raison>1}
\sum_{m\in \Z^*}e^{-(\frac{1}{2}+\epsilon)d_{\H^N}({\bf o}, p^{m}\cdot {\bf 
o})}
\sum_{n\in \Z^*}e^{-(\frac{1}{2}+\epsilon)d_{\H^N}({\bf o}, h^n\cdot {\bf 
o})}>1,
\end{equation} 
which proves that the critical exponent of 
$G$ is strictly greater than ${1\over 2}$ for $a=+\infty$.

The same property holds in fact for finite but large enough values of $a$. Indeed,  there exists $m_a\geq 1$, with $m_a\to +\infty$ as $a \to +\infty$, such that the geodesic segments $[{\bf o}, p^m\cdot {\bf o}]$ for $-m_a\leq m\leq m_a$ remain inside the stripe $\R^{N-1}\times [0, a] \subset \cH$ corresponding to the cuspidal end of type $(a, \alpha)$ 
; since $g_a$ has curvature   $-1$ in this stripe (see Remark \ref{stripe}), the quantities $d_a ({\bf o}, p^m\cdot {\bf o})-2\ln \vert m\vert $ remain also bounded  for these values of $m$, uniformly in $a$ (\cite{DOP}, lemme  4). In the same way, for any $a\geq 0$ the geodesic segments $[{\bf o}, h^{n}\cdot {\bf o}]$  lie  in the area of constant curvature of the metric $g_a$ so that $d_{a}({\bf o}, h^{n}\cdot {\bf 
o})= d_{\H^N}({\bf o}, h^{n}\cdot {\bf 
o})$.
So, by (\ref{raison>1}), for
$\epsilon \in ]0, \epsilon_0]$, one gets 
$$\liminf _{a\to +\infty}
\sum_{|m|\leq m_{a}}e^{-(\frac{1}{2}+\epsilon)d_{a}({\bf o}, 
p^m\cdot {\bf o})}
\sum_{n\in \Z^*}e^{-(\frac{1}{2}+\epsilon)d_{a}({\bf o}, h^{n}\cdot {\bf 
o})} >1.$$
There thus exists $a_0>0$ such that, for $a\geq a_0$ one gets  
$\ds  \sum_{g\in G}\e^{-(\frac{1}{2}+\epsilon) d_a({\bf o}, g\cdot {\bf o})}=+\infty.$
 This last inequality implies that  $\delta_{G}>{1\over 2}$ when  $a \geq a_0$ ; by (\cite{DOP}, Proposition 1), the group $G$ is thus divergent since $\delta_{\langle p\rangle}={1\over 2}$.\fdem

We now want to check that  there exists  some  $a \in ]0, 
a_{0}[$ 
such that  the group $G$ 
 is divergent with 
 $\delta_{G}= {1\over 2}$ when $X$ is endowed with the metric $g _a$ ; to prove this,  we need to compare  the Poincar\'e
 series ${\bf P}_{G}(s)$  with  the potential of some Ruelle 
 operator ${\mathcal L}_{a, s}$  associated with $G$ that we introduce in the following paragraph.
 
 {\bf From now on, we fix a Schottky group $G=\langle p, h \rangle$ satisfying the conclusions of  Proposition  \ref{Schottkygroupsconvergentdivergent} and subsets ${\mathcal U}_p$ and ${\mathcal U}_h$ in $X\cup \partial X$ satisfying   conditions $(1), (2), (3)$ and $(4)$ above.}
 \subsection{Spectral radius of the Ruelle  operator and Poincar\'e exponent}
  We   introduce the  family  $({\mathcal L}_{a, s})_{(a, s)} $   of  {\it Ruelle operators}  associated with $G=\langle p, h \rangle$ defined  formally by :
for any $a \in [0, a_0], s \geq 0,  x \in \partial X$ and  any bounded Borel function $\phi: \partial X \to \R$
\begin{equation}\label{defiruelle}
 {\mathcal L}_{a, s}\phi(x) = \sum_{\gamma \in \{p, h\}}\sum_{n \in 
 \Z^*}1_{x \notin U_{\gamma}} e^{ -s {\mathcal 
 B}_{x}^{(a)}(\gamma^{-n}\cdot {\bf o},{\bf o})}\phi(\gamma^n\cdot x).
 \end{equation}
 The sequence $(p^n\cdot {\bf o})_{n}$ accumulates at $\xi_{p}$. 
So, for 
 any $x\in {\mathcal U}_{h}$ the angle at  ${\bf o}$ of the triangle 
 $\widehat{{\bf x}\ {\bf o}\  p^n\cdot {\bf o}}$ 
 is greater than $\theta/2$ for $n$ large enough and  the sequence 
 $({\mathcal  B}_{x}^{(a)}(p^{-n}\cdot{\bf o}, {\bf o})-d_{a}({\bf o}, 
p^n\cdot {\bf o}))_{n}$ is bounded uniformly 
 in $x \in {\mathcal U}_{h}$ and $a \geq 0$. In the same way, the 
sequence
 $({\mathcal  B}_{x}^{(a)}(h^{-n}\cdot {\bf o}, {\bf o})-d_{a}({\bf o}, 
h^n\cdot {\bf o}))_{n}$ is bounded uniformly 
 in $x \in {\mathcal U}_{p}$ and $a \geq 0$. It readily implies that
 ${\mathcal L}_{a, s}\phi(x)$ is finite when $s \geq \max( \delta_{\langle p\rangle}, \delta_{\langle h\rangle})= {1\over 2}$ 
and 
 that it acts on the space  ${\mathbb L}^{\infty}(\partial X)$ of 
bounded  Borel functions on $\partial X$.

By a similar argument, for any $k\geq 1$, the 
quantities
$$ \Bigl({\mathcal  B}_{y}^{(a)}((p^{m_{1}}h^{n_{1}}\cdots p^{m_{k}}
h^{n_{k}})^{-1}\cdot {\bf o}, {\bf o})
-d_{a}({\bf o}, p^{m_{1}}h^{n_{1}}\cdots p^{m_{k}}h^{n_{k}}\cdot {\bf o})\Bigr)_{n}$$
and  
$$\Bigl({\mathcal  B}_{x}^{(a)}((h^{n_{1}}p^{m_{1}}\cdots h^{n_{k}} 
p^{m_{k}})^{-1}\cdot {\bf o}, {\bf o})
-d_{a}({\bf o}, h^{n_{1}}p^{m_{1}}\cdots h^{n_{k}}p^{m_{k}}\cdot {\bf o})\Bigr)_{n}$$\ 
are  bounded 
uniformly in $m_{1}, n_{1}, \cdots, m_{k}, n_{k} \in \Z^*$ and 
$x \in {\mathcal U}_{p}, y \in {\mathcal U}_{h}.$ One thus gets
$$
\sum_{m_{1}, \cdots n_{k}\in \Z^*} \exp(d_{a}({\bf o}, 
p^{m_{1}}h^{n_{1}}\cdots 
p^{m_{k}} h^{n_{k}}\cdot {\bf o}))
\asymp
|{\mathcal L}_{a,s}^{2k}1|_{\infty}$$
 which states that the series ${\bf P}_{G}(s)$ and 
$\displaystyle{\sum_{k \geq 1}|{\mathcal L}_{a,s}^{2k}1|_{\infty}}$ 
diverge or converge simultaneously.
Now, the fact that ${\mathcal L}_{a,s}$ is a non negative operator implies that  the limit 
$\displaystyle{\lim_{k \to +\infty}\Bigl(|{\mathcal 
L}_{a,s}^{2k}1|_{\infty}\Bigr)^{1\over 2k}}$ is equal to the spectral 
radius 
$\rho_{\infty}({\mathcal L}_{a,s})$ 
of ${\mathcal L}_{a,s}$ on ${\mathbb L}^{\infty}(\partial X)$.
 We have  thus 
established  the

\begin{fact}\label{potentielRuelle}
The Poincar\'e series ${\bf P}_{G}(s)$
 and the potential
 $\displaystyle{\sum_{k \geq 1}|{\mathcal L}_{a,s}^{k}1|_{\infty}}$ 
diverge or converge simultaneously. In particular,  if $\delta_a$ denotes the Poincar\'e
exponent 
  of $G$  for the metric $g_{a}$, one gets  
  $$
  \delta_a=\sup\Bigl\{s \geq 0 : \rho_\infty (\mathcal L_{a, s})\geq 1\Bigr\}=\inf 
  \Bigl\{s \geq 0 : \rho_\infty (\mathcal L_{a, s})\leq 1\Bigr\}.
  $$
\end{fact}
Consequently, since $G$ satisfies the conclusions of Proposition
 \ref{Schottkygroupsconvergentdivergent},  one gets 
 \begin{itemize}
 \item the series
 $\displaystyle{\sum_{k \geq 1}|{\mathcal L}_{0,{1/2}}^{2k}1|_{\infty}}$  converges ;
 \item   for  $a_0 $ large enough, the series $\displaystyle{\sum_{k \geq 1}|{\mathcal L}_{a_0,{1/2}}^{2k}1|_{\infty}}$ diverges
 \end{itemize}
 which implies in particular $\rho_\infty({\mathcal L}_{0,{1/2}}) \leq 1$  and $\rho_\infty({\mathcal L}_{a_0,{1/2}}) \geq 1$. 
 
 We will prove that for some value  $a_* \in ]0, a_0[$ one gets  $\rho_\infty({\mathcal L}_{a_*,{1/2}})= 1$ ; the unicity of $a^*$ will be specified in the last section.

We  first need to control the regularity of the function $a \mapsto \rho_\infty({\mathcal L}_{a_0,{1/2}})$. It will be quite simple to check that the  function $a \mapsto {\mathcal L}_{0,{1/2}}$ 
is continuous from $\R^+$ to the space of bounded operators on 
${\mathbb L}^{\infty}(\partial X)$ ; unfortunately, the function $\mathcal L\mapsto \rho_\infty(\mathcal L)$ is in general only lower semi-continuous. In the case of the family of Ruelle operators we consider here,  this function  will be in fact continuous,  because of the very special form of the spectrum in this situation.

\subsection{On the spectrum of the Ruelle operators}

Throughout this section,  we will use the following 
\begin{notation} For any $a \in [0, a_0], x \in \partial X$ and $\gamma \in G$,  we will set 
\begin{itemize}
\item 
${\mathcal L}_{a}= {\mathcal L}_{a, 1/2}$ 
and 
$ \rho_{\infty}(a) = \rho_{\infty}({\mathcal L}_{a, 1/2})$.
\item  $b_{a}(\gamma,x)={\mathcal B}_{x}^{(a)}(\gamma^{-1}{\bf o}, 
{\bf 
o})$
\end{itemize}

Furthermore, $\delta_a$ will denote the Poincar\'e exponent of $G$ with respect to the metric $g_a$ and, for $\gamma \in \{p, h\}$ and $n \in \Z^*$,   the ``weight'' function $w_{a}(\gamma^n, .)$ is  defined by 
\begin{eqnarray*}
w_{a}(\gamma^n, .) \quad : \quad  \partial X & \to & \R^+\\
 x &\mapsto&
 1_{x \notin U_{\gamma}}e^{-\delta_{a}b_{a}(\gamma^n,x)}
\end{eqnarray*}
\end{notation}

With these notations, the Ruelle operator $\mathcal L_a$ introduced in the previous paragraph may be expressed as follows : for any $\phi \in {\mathbb L}^\infty(\partial X)$ and any $x \in \partial X$, 
\begin{equation}\label{defiruellebis}
 {\mathcal L}_{a}\phi(x) = \sum_{\gamma \in \{p, h\}}\sum_{n \in 
 \Z^*}w_{a}(\gamma^n, x)\phi(\gamma^n\cdot x).
 \end{equation}
The iterates of $\mathcal L _a$ are given by 
\begin{equation}\label{defiruelle iterates}{\mathcal L}_{a}^k\phi(x) = \sum_{\gamma \in 
G(k)}w_{a}(\gamma,x)\phi(\gamma\cdot x)
\end{equation}
where $G(k)$  is the set of $\gamma = \alpha_{1}^{n_{1}}\cdots 
 \alpha_{k}^{n_{k}} \in G$  of length  $k$ (with $ \alpha_{i+1}\neq 
 \alpha_{i}^{\pm 1}$) and $w_{a}(\gamma,x)=1_{x\notin U_{\alpha_k}}e^{-\delta_{G}b_a(\gamma,x)} $  when $\gamma$ has last letter $\alpha_k$ ; observe that  we have the following ``multiplicative  cocycle property'' :
\begin{equation}\label{mult-cocycle}
w_{a}(\gamma,x)= \prod_{i=1}^k w_{a}(\alpha_{i}^{n_{i}}, 
\alpha_{i+1}^{n_{i+1}}\cdots \alpha_{k}^{n_{k}}\cdot x).
\end{equation}

 We  will see that the 
 ${\mathcal L}_{a}, a \geq 0$, act on the space 
 $C(\partial X)$ of real  valued  continuous functions on 
 $\partial X$ and that the map $a\mapsto \mathcal L_a$ is continuous.  Nevertheless,    the function 
$a \mapsto \rho_{\infty}(a)$ is only lower semi-continuous in general and may present   discontinuities. The main idea to avoid this difficulty  is to introduce a 
Banach 
space  on which the ${\mathcal L}_{a}$ act  quasi-compactly. 

In the sequel, we will  consider the  restriction of 
the ${\mathcal L}_{a}$ to some  subspace
of $C(\partial X)$ of H\"{o}lder continuous functions.

\begin{notation}   We denote  $\mathbb L_{a, \omega}(\partial X)$
  the space of   H\"{o}lder continuous functions on $\partial X$ defined by
$$
\mathbb L_{a, \omega}(\partial X) := \{\phi \in C(\partial X) / 
 |\phi |_{a, \omega}= |\phi |_{\infty}+[\phi]_{a, \omega}<+\infty\}
 $$
 where
  $\ds{[\phi]_{a, \omega}= \sup_{\gamma\in \{p,h\}}
 \sup_{\stackrel{x, y \in U_{\gamma}}{x\neq y}}
 \frac{|\phi(x)-\phi(y)|}{D_a(x,y)^{\omega}}}$ denotes the 
$ \omega$-H\"{o}lder
coefficient of $\phi$ with respect to the distance $D_a$.

When $a=0$ we will omit the index $D_0$ and set 
$ {\mathbb L}_{\omega}(\partial X) := {\mathbb L}_{ 0,\omega}(\partial X) $.
\end{notation}
The spaces $(\mathbb L_{a, \omega}(\partial X),  |. |)$ are 
$\C$-Banach space  and the identity map  
from $(\mathbb  L_{a, \omega}(\partial X),  |. |_{a, \omega})$ to $(C(\partial X), 
 |. |_{\infty})$ are compact. 
 
 We now want to prove that each operator $\mathcal L_a$ acts on $\mathbb L_{a, \omega}(\partial X) $; in fact,  we need a stronger result, i.e that each $\mathcal L_a$, for $0\leq a\leq a_0$, acts on  ${\mathbb L}_{\omega}(\partial X) $. It will be a direct consequence  of the following :

\begin{lemma}\label{poidsequiholder}  There exists  $\omega_0 \in ]0, 1[$ such that 
for any $\omega \in ]0, \omega_0],$ any $\gamma \in \{p,h\}$, any $a \in [1, a_{0}]$ and any $n \in \Z^*$, the 
function 
$w_{a}(\gamma^n, .)$ belongs to $\mathbb L_\omega(\partial X)$; furthermore, 
the 
sequence 
$\bigl(e^{\delta_{G} d_{a}({\bf o}, \gamma^n \cdot {\bf o})}
 |w_{a}(\gamma^n, 
.) |_{\omega}\bigr)_{n}$
is bounded.
\end{lemma}
 \noindent Proof.  The cluster points of the sequence
 $(\gamma^n \cdot {\bf o})_{n}$ belong to $U_{\gamma}$. Since the curvature  
 is pinched, the quantity $b_{a}(\gamma^{n}, x))-d({\bf 
 o}, \gamma^{n} \cdot {\bf o} )$ is  bounded uniformly in $n \in \Z^*, x 
\in 
 \partial X\setminus U_{\gamma}$ and $a \in [0, a_{0}]$ ; so is the sequence 
 $\bigl(e^{\delta_{G} d_{a}({\bf o}, \gamma^n \cdot {\bf o})} 
|w_{a}(\gamma^n, 
.) |_{\infty}\bigr)_{n}$. In order to 
 control the $\omega$-H\" older coefficient  of  $w_{a}(\gamma^n, 
.)$, we use Fact \ref{Fact} : the   functions 
$x \mapsto {\mathcal  B}_{x}^{(a)}({\bf o}, \gamma^{-n} \cdot {\bf o})$  are
equi-Lipschitz continuous on $\partial X\setminus U_{\alpha}$ with 
respect to $D_{a}$, since 
once again the cluster points of the sequence
 $(\gamma^n \cdot {\bf o})_{n}$ belong to $U_{\gamma}$. 
 More precisely, the sequence $\bigl(e^{\delta_{G} d_{a}({\bf o}, \gamma^n \cdot {\bf o})}
 |w_{a}(\gamma^n, 
.) |_{\omega}\bigr)_{n}$
is bounded for any $a \in [0, a_0]$ and 
 for $\omega$ given  by inequality  (\ref{equi-holder}) 
.\fdem

{\bf From now and for once , we fix $\omega_0 \in ]0, 1[$ satisfying the conclusion of the above Lemma.} We know that,  for $a \in [0, a_0]$, the operator  
 ${\mathcal L}_{a}$ acts on 
$\mathbb L_\omega(\partial X)$  whenever  $\omega \in ]0, \omega_0]$ ; let  $\rho_{\omega}(a)$ denote  the 
spectral radius of ${\mathcal L}_{a}$ on   $\mathbb L_\omega(\partial X)$.
We have the 
\begin{proposition}
 For any $\omega \in ]0, \omega_0]$ and $a \in [1, a_{0}]$, 
one gets 
\begin{itemize}
\item 
$\rho_{\omega}(a)=\rho_\infty(a)$
\item 
$\rho_{\omega}(a)$ 
is a simple   eigenvalue of the operator  ${\mathcal L}_{a}$  
acting on $\mathbb  L_{\omega}(\partial X)$ and  the 
associated eigenfunction is non negative  on  $\partial X$.
\end{itemize}

 Furthermore, the operator
 ${\mathcal L}_{a}$  is quasi-compact on 
$\mathbb  L_{\omega}(\partial X)$ : 
there exists $r<1$ such that the essential spectral  radius of 
${\mathcal L}_{a}$  on 
$\mathbb  L_{\omega}(\partial X)$ is less than $r   \rho_{\omega}(a)$. 

In 
particular the eigenvalue  $\rho_{\omega}(a)$ is  isolated in 
the spectrum of ${\mathcal L}_{a}$,  it is simple  and the corresponding   eigenfunction
is non-negative.
\end{proposition}
 \noindent Proof.
 Fix   $x, y \in \partial X\cap 
U_{p}$ ; one gets
\begin{eqnarray*}
|{\mathcal L}_{a}^k\phi(x) -{\mathcal L}_{a}^k\phi(y)|
&\leq& \sum_{\gamma \in 
G(k)}w_{a}(\gamma,x)|\phi(\gamma\cdot x)-\phi(\gamma\cdot y)|\\
& & \qquad \qquad +
\sum_{\gamma \in 
G(k)}|w_{a}(\gamma,x)-w_{a}(\gamma,y)|\times  |\phi |_{\infty}. 
\end{eqnarray*}
Note that in these sums, it is sufficient to consider the 
$\gamma \in G(k)$ with last letter  $\alpha_{k}\neq p$. 
For such $\gamma$ the quantity 
$b_{a}(\gamma, x)$ is greater than
$d_{a}({\bf o}, \gamma\cdot {\bf o})-C$ for some constant $C$  which 
depends only on the angle $\theta_{0}$ and the bounds on the 
curvature ; in particular $b_{a}(\gamma, x)\geq 1$ for all but 
finitely many $\gamma \in G$ with last letter $\neq p$. It 
readily follows that 
$\ds{\liminf_{\stackrel{\gamma \in G(k)}{k\to +\infty}}
\frac{b_{a}(\gamma,x)}{k}>0}$, uniformly in $x\in U_{p}$. 
In other words, there exists $0<r<1$ and $C>0$ such that 
$$|\gamma'(x)|_{a}\leq Cr^k$$
 for any $ k \geq 1, x \in U_{p}$ and
$\gamma \in G(k)$ 
with last letter $\neq p$.  
The same argument works when $x, y \in \partial X\cap U_h$.

We thus obtain  the inequality
$$
[{\mathcal L}_{a}^k\phi]_{\omega}\leq 
r_{k}[\phi]_{\omega}+R_{k}|\phi|_{\infty}
$$
with $r_{k} = C K_{\omega}^2 r^k |{\mathcal L}_{a}^k|_{\infty}$ 
and $R_{k} = 
\sum_{\gamma \in G(k)}[w_{a}(\gamma, .)]_{\omega}$.

Note that  ${\mathcal 
L}_{a}$ is a non-negative operator, so that the quantity 
$\limsup_{k}|{\mathcal 
L}_{a}^k1|_{\infty}^{1/k}$ is equal to  the spectral radius 
$\rho_{\infty}(a)$ of ${\mathcal 
L}_{a}$ on $C(\Lambda)$. Using a  version due to H. Hennion  of the 
Ionescu-Tulcea-Marinescu's  theorem concerning quasi-compact 
operators, one may conclude that the essential spectral radius of 
${\mathcal 
L}_{a}$ on $\mathbb  L_{\omega}(\partial X)$ is less than $r 
\rho_{\infty}(a)$ ; in other words, the spectral values of 
${\mathcal 
L}_{a}$ with modulus $\geq r\rho_{\infty}(a)$ are isolated 
eigenvalues with finite multiplicity in the spectrum of 
${\mathcal 
L}_{a}$.
 This implies
 in particular that $\rho_{\omega}(a) =  \rho_{\infty}(a)$. The 
inequality
 $\rho_{\omega}(a) \geq   \rho_{\infty}(a)$ is obvious since the 
function $1$ 
 belongs to ${\mathbb L}_{\omega}(\partial X)$.
Furthermore, the strict inequality  
would imply   the existence of  a function $\phi \in  {\mathcal 
L}_{a}$ such that ${\mathcal 
L}_{a}\phi = \lambda \phi$ for some $\lambda \in \C$ of modulus 
$>  \rho_{\infty}(a)$ ; this would give $| \lambda| |\phi| \leq
{\mathcal 
L}_{a} |\phi|$ so that $|\lambda |\leq \rho_{\infty}(a)$, which is 
a contradiction. 

It remains to  control the value $\rho_{\omega}(a)$ in 
the spectrum of ${\mathcal L}_{a}$.  The operator ${\mathcal 
L}_{a}$ being  non negative and compact on $\mathbb L_\omega(\partial X)$, its spectral radius   $\rho_{\omega}(a)$ is an  eigenvalue with  associated 
eigenfunction 
$\phi_{a} \geq 0$.

 Assume that $\phi_{a}$ vanishes at $x_{0} \in 
\partial X$ and let $g \in \{p, h\}$  such that
$x_{0}\in U_{g}$; 
the equality ${\mathcal L}_{a}\phi_{a}(x_{0})= 
\rho_{\omega}(a)\phi_{a}(x_{0})$ implies 
that  $\phi_{a}(\gamma \cdot x_{0}) = 0$ for any $\gamma \in G$ with 
last letter  $\neq g$.
 By minimality of the action of $G$ 
on $\partial X$ one  thus has $\phi_{a}= 0$ on 
$\partial X$. Consequently, the function $\phi_{a}$ is  non negative.

Let us now check that $\rho_{\omega}(a)$ is a simple eigenvalue of $\mathcal L_a$. Consider the  operator   $P$  
defined formaly by  $\ds{P( f ) = 
\frac{1}{\rho_{\omega}(a)\phi_{a}}{\mathcal 
L}_{a}(f \phi_{a})}$ ; this operator is well defined  on $\mathbb L_\omega(\partial X)$ since 
$\phi_{a}$ does not vanish, it  is non negative, 
quasi-compact with   spectral radius 
 $1$ and  Markovian (that is to say $P1 = 1$).  If  $f \in {\mathbb L}_{\omega}(\partial X)$
satisfies the equality $Pf=f$ one  considers a point $x_{0}
\in \partial X$ such that $|f(x_{0})|= |f |_{\infty}$ and $g \in \{p, h\}$ such  that  $x_{0}\in U_{g}$. 
An argument of convexity applied to 
the inequality  $P|f|\leq |f|$ readily implies 
$|f(x_{0})|= |f(\gamma\cdot x_{0})|$ for any $\gamma \in G$ 
with last letter $\neq g$ ; by minimality of the action of
$G$ on $\partial X$ it follows that  the modulus of $f$ is
constant  
on $\partial X$. Applying  again an argument of convexity and 
the
minimality of the action of
$G$ on its limit set, one 
proves that $f$ is  in fact constant on $\partial X$ ; it follows that $\mathbb C \phi_a$ is the eigenspace associated with $\rho_\omega(a)$ on $\mathbb L_\omega(\partial X)$.\fdem

\subsection{Regularity of the function $a \mapsto {\mathcal L}_{a}$}

In this section we will establish the following

\begin{proposition}
 For any $0< \omega <  \omega_{0}$,  the function 
$a \mapsto {\mathcal L}_{a}$ is continuous from $[1, 
a_{0}]$ to the space 
of continuous linear operators on  
$(\mathbb L_{\omega}(\partial X), |. |_{\omega})$.
\end{proposition}
 \noindent Proof. It suffices to check that,  for $\gamma \in 
\{h,p\}, 
a, a' \in [1,a_{0}]$ and $0< \omega <  \omega_{0}$ one has 
$$\lim_{a'\to a}\sup_{n \in \Z}e^{\delta_{G}d_{0}({\bf o}, 
\gamma^n \cdot {\bf o})} |w_{a'}(\gamma^n,.)-w_{a}(\gamma^n,.) |_{\omega}=0.$$
First one gets 
\begin{eqnarray*}
    |w_{a'}(\gamma^n, .)-w_{a}(\gamma^n, .)|
    &=& e^{-\delta_{G}b_{a}(\gamma^n, .)}
    |e^{-\delta_{G}(b_{a'}(\gamma^n,.)-b_{a}(\gamma^n,.))}-1|\\
    &\leq&
    C e^{-\delta_{G}d_{a}(\gamma^n{\bf o}, {\bf o})}
     |e^{-\delta_{G}(b_{a'}(\gamma^n,.)-b_{a}(\gamma^n,.))}-1|  
\end{eqnarray*}
where the constant $C$ depends only one the bounds on
the curvature. 

Since the axis 
of $h$ lies in the area of $X$ where the curvature is $-1$, the quantity
$d_{a}({\bf o}, h^n\cdot  {\bf o})-|n|l_{h}$, where $l_{h}$ denotes the hyperbolic
lenght of the closed geodesic associated with $h$, is bounded uniformly in $a \in [0, a_0]$ and $n \in \Z^*$ ; the same holds for 
the quantity
$d_{a}(p^n\cdot {\bf o}, {\bf o})-d_{0}(p^n\cdot {\bf o}, {\bf o})$. Consequently
$$
|w_{a'}(\gamma^n, .)-w_{a}(\gamma^n, .)|\leq
C' e^{-\delta_{G}d_{0}(\gamma^n{\bf o}, {\bf o})}
|e^{-\delta_{G}(b_{a'}(\gamma^n,.)-b_{a}(\gamma^n,.))}-1|$$
and we have to study the regularity  of the function $a\mapsto 
b_{a}(\gamma^n,x)$,
for any point $x \notin U_{\gamma}$.
By inequalities (\ref{gromovquasi-iso}), one gets
$$(y|z)_{a'}\to (y|z)_{a}\quad \mbox{\rm as} \quad a'\to a,$$
when $
(y|z)_{a}
$ remains  bounded. There are thus two cases to consider: 
\begin{itemize}
\item We first consider  the case $\gamma = p$. For any $n \in \Z^*$ 
  let $y_{n}$ be the point in $\partial X$ such that  ${\bf o}$
  belongs to  the geodesic 
  ray $[p^n\cdot {\bf o}, y_{n})$ (for the metric $g_{a}$) ; 
this ray is in fact  a quasi-geodesic for any $a' \in [0, a_{0}]$,  
so  the 
point  ${\bf o}$ belongs to some bounded 
neigbourhood of the geodesic ray (for $g_{a'}$)  from $p^n\cdot {\bf o}$ to 
$x_{n}$ (
note that  $\ds{\inf_{n\in \Z^*}
D_{0}(y_{n}, \xi_{p})>0}$ by convexity of the horospheres).
For any $n \in \Z^*$ and $a' \in [0, a_{0}]$  one gets
$$b_{a'}(p^n,x) = 
(p^n\cdot x|y_{n})_{a'}-(x|p^{-n}\cdot y_{n})_{a'}-b_{a'}(p^n,p^{-n}\cdot y_{n}).$$
Since $p^n\cdot x \to \xi_p$ as $|n|\to+\infty$  and $\ds{\inf_{n\in \Z^*}
D_{0}(y_{n}, \xi_p)>0}$, one gets
$$(p^n\cdot x|y_{n})_{a'}\to (p^n\cdot x|y_{n})_{a}$$ as $a'\to a$, uniformly in $n 
\in \Z^*$ and $x \notin {\mathcal U}_{p}$. In the same way,
since $\ds{\inf_{n\in \Z^*}
D_{0}(y_{n}, \xi_p)>0}$, the sequence 
$(p^{-n}\cdot y_{n})_{n}$ converges to $\xi_p$ as $|n| \to +\infty$
so that 
$(x|p^{-n}\cdot y_{n})_{a'}
\to (x|p^{-n}\cdot y_{n})_{a}$ uniformly in $n 
\in \Z^*$ and $x \notin {\mathcal U}_{p}$. Atlast one has
$b_{a'}(p^n,p^{-n}\cdot x_{n})={\mathcal B}^{(a')}_{x_{n}}({\bf o}, p^n\cdot {\bf 
o})
= d_{a'}({\bf o}, p^n\cdot {\bf o})$ ; the geodesic segment
$[{\bf o}, p^n\cdot {\bf o}]$ is included in the horosphere
${\mathcal H}$, so that
$$b_{a'}(p^n,p^{-n}\cdot x_{n}) \to
b_{a}(p^n,p^{-n}\cdot x_{n})$$ as $a'\to a$, uniformly in $n \in \N^*.$ 
\item Consider now  the case when $\gamma = h$  ; for any $n \geq 1, $ one gets
$$b_{a}(h^n,x) = 
(h^n\cdot x|h^n\cdot \xi_{h}^+)_{a}-(x|\xi_{h}^+)_{a}-b_{a}(h^n,\xi_{h}^+)$$
with
$b_{a}(h^n,\xi_{h}^+)=nl_{h}$.
The facts that $x\notin U_{h}$ and $\xi_{h}^+\in U_{h}$ readily implies
$(x|\xi_{h}^+)_{a'} \to (x|\xi_{h}^+)_{a}$ as $a'\to a$.
On the other hand $h^n\cdot x\to x_{+}$ as $n\to 
+\infty$
    so that  $(h^n\cdot x|h^n\cdot \xi_{h}^+)_{a}\to (x_{+}|\xi_{h}^+)_{a}$ ; since
  $\xi_{h}^+\neq x_{+}$, the Gromov product  $(x_{+}|\xi_{h}^+)_{a}$ is equal 
to 
  $-\log d_{a}({\bf o}, (x_{+}\xi_{h}^+))$  up to a bounded term and 
  the sequence $((h^n\cdot x|h^n\cdot \xi_{h}^+)_{a})_{n\geq 1}$ 
is  bounded uniformly in 
  $a\in [0,a_{0}], x\notin U_{h}$ and $n \in \N$.
  It 
  readily follows that $(h^n\cdot x|h^n\cdot\xi_{h}^+)_{a'}\to (h^n\cdot x|h^n\cdot\xi_{h}^+)_{a}$ 
  as $a'\to a$, for any $n \geq 1$. A similar argument holds for $n 
  \leq -1$. Finally, $b_{a'}(h^n,x)-b_{a}(h^n,x)\to 0$ uniformly in 
  $n \geq 0$ and $x \notin U_{h}$ and the lemma is proved for $\gamma 
  =h$.
\end{itemize}
Finally one has proved that
for $\gamma \in \{h,p\}$ and $ 
a, a' \in [0,a_{0}]$  one has 
$$\lim_{a'\to a}\sup_{n \in \Z}e^{\delta_{G}d_{0}({\bf o}, 
\gamma^n \cdot {\bf o})} |w_{a'}(\gamma^n,.)-w_{a}(\gamma^n,.) |_{\infty}=0.$$
To achieve the \noindent Proof., we use the classical fact that if a bounded sequence $(f_n)_n$  in $\mathbb L_{\omega_0}(\partial X)$ converges uniformly to some  (continuous) function $f$, then the convergences remains valid in
 $\mathbb L_{\omega}(\partial X)$  for any $0<\omega<\omega_0$ : 
namely, we may fix  $\epsilon 
>0$ and note that, for $0<   \omega \leq  \omega_{0}$, the following 
inequality holds
$$
[w_{a'}(\gamma^n,.)-w_{a}(\gamma^n,.)]_{\omega}\leq 
\frac{2 |w_{a'}(\gamma^n,.)-w_{a}(\gamma^n,.) 
|_{\infty}}{\epsilon^ \omega}
+ 
[w_{a'}(\gamma^n,.)-w_{a}(\gamma^n,.)]_{ \omega_{0}}\epsilon^{ \omega_{0}- \omega} 
$$
which immediately gives 
$$
 |w_{a'}(\gamma^n,.)-w_{a}(\gamma^n,.) |_{\omega}\leq 
\bigl(\frac{2}{\epsilon^{ \omega}}+1 \bigr)
 |w_{a'}(\gamma^n,.)-w_{a}(\gamma^n,.) |_{\infty}
+ |w_{a'}(\gamma^n,.)-w_{a}(\gamma^n,.) |_{ \omega_{0}}
\epsilon^{ \omega_{0}- \omega}. 
$$
One achieves the \noindent Proof. letting $a' \to a$ and $\epsilon \to 0$.\fdem

\subsection{\noindent Proof. of the Main Theorem }
We are now able to achieve the \noindent Proof. of the Main Theorem.   We fix $ \omega 
\in ]0,  \omega_{0}[.$ 

Since the 
spectral radius $\rho_{\omega}(a)$ 
of the operator ${\mathcal L}_{a}$  acting on ${\mathbb L}_{\omega}$ 
is  an eigenvalue and is 
isolated is the spectrum of  ${\mathcal L}_{a}$, the function $a \mapsto 
\rho_{\omega}(a)$
has the same regularity than   $a \mapsto {\mathcal L}_{a}$ ; it is  thus 
continuous on $[1, a_{0}]$.
Furthermore, for any $a \in [1, a_{0}]$,
the eigenfunction $\phi_{a}$ associated with 
$\rho_{\omega}(a)$ is non negative on $\partial X$. So one has  
$\phi_{a}\asymp 1$,  which readily implies that 
$
 \displaystyle{\vert {\mathcal L}_{a}^{2k}\phi_{a}\vert _{\infty}
 \asymp
 \vert {\mathcal L}_{a}^{2k}1\vert_{\infty}}$ uniformly in $k\geq 1$.
By the equality ${\mathcal L}_{a}\phi_{a} = \rho_{\omega}(a)\phi_{a}$, it 
follows that $\rho_{\omega}(a)=\rho_{\infty}(a)$. 

 By the choice of the metrics $g_{a}$, we have  
$ \rho_{\infty}(0) \leq 1$ 
and 
$ \rho_{\infty}(a_{0}) \geq 1$ ; so there exists
$a_{*}\in ]0, a_{0}[$ such that 
$\rho_{\omega}(a_{*}) =\rho_{\infty}(a_{*})  = 1$.  

On the other hand, the function 
$s  \mapsto \rho_{\omega}({\mathcal L}_{a_{*},s})$ is strictly decreasing on 
$\R^{+}$. Fix  $s>\delta_{\langle p\rangle}$ ; one has
$\rho_{\omega}(a_{*})<1$ and   the series
$P_{G}(s)$  thus  converges  when $X$ is endowed with the 
metric 
$g_{a^{*}}$. This proves that for the value $a_{*}$ of the parameter 
$a$ 
the critical exponent of $G$ is less than  $\delta_{\langle p\rangle} $; since $p 
\in G$,
one has in 
fact $\delta_{G}=\delta_{\langle p\rangle}$.

Atlast, since $\phi_{a_{*}}\asymp 1$, 
 one has  
 $\displaystyle{\sum_{k \geq 1}|{\mathcal 
L}_{a_{*},\delta_{\langle p\rangle}}^{2k}1|_{\infty}
 \asymp  \sum_{k \geq 1}|{\mathcal 
L}_{a_{*},\delta_{\langle p\rangle}}^{2k}\phi_{a_{*}}|_{\infty}}$
 ;  these two series diverge  in fact because of 
 the equality
 ${\mathcal L}_{a_{*},\delta_{\langle p\rangle}}\phi_{a_{*}}=\phi_{a_{*}}$.
By the Fact \ref{potentielRuelle}, it follows that  for the value $a_{*}$ of the 
parameter $a$, the series $P_{G}(\delta_{G})$ diverges.

By criteria (\ref{BMfinie}), one easily sees that $m_\Gamma$ is finite when $\alpha >2$ and infinite when $\alpha \in ]1, 2]$.

This achieves 
the \noindent Proof. of the Main Theorem.
\fdem

\subsection{Complement}
A natural question  is the one of unicity of the value $a_*$ of the parameter $a$ such that the spectral radius $\rho_\infty(a)$ of $\mathcal L_a$ is equal to $1$; this unicity is  not necessary to prove the main Theorem but nevertheless it is of interest to describe for instance the behavior of the orbital function of $G$ when $a$ varies. It will be the subject of a forecoming work.

By the continuity of the function $a\mapsto \rho_\infty(a)$, the unicity of $a_*$ is a direct consequence of the  strict monotonicity of this function. 
We thus have to prove that $\rho({\mathcal L}_{a})<\rho({\mathcal L}_{a'})$ for any $a, a'$ in $[0, a_0]$ such that $a<a'$. Note first that, for any fixed ${\bf x} \in X$ one gets 
$$
\rho({\mathcal L}_a)=\rho_{\infty}({\mathcal L}_a)=\lim_{k\to +\infty} \Bigl(\Big\Vert
\sum_{\stackrel {\gamma \in \Gamma_{2k}}{l(\gamma)=h}} e^{-{1\over 2}{\mathcal B}^{(a)}_{.}(\gamma ^{-1}\cdot{\bf x},  \bf x)}
\Big\Vert_\infty\Bigr)^{1\over 2k}= \lim_{k\to +\infty} \Bigl(
\sum_{\stackrel {\gamma \in \Gamma_{2k}}{l(\gamma)=h}} e^{-{1\over 2}d_a({\bf x}, \gamma\cdot{\bf x})}
\Bigr)^{1\over 2k}
$$
and we have thus to check  that there exists $C>0$ and $\rho:=\rho(a, a')<1$ such that, for any $n \geq 1$, one gets
\begin{equation}\label{ineq-stricte}
\sum_{\stackrel {\gamma \in \Gamma_{2k}}{l(\gamma)=h}} e^{-{1\over 2}d_a({\bf x}, \gamma\cdot{\bf x})}    
\leq  C \rho^k \sum_{\stackrel {\gamma \in \Gamma_{2k}}{l(\gamma)=h}} e^{-{1\over 2}d_{a'}({\bf x}, \gamma\cdot{\bf x})}.
\end{equation}
For any $x \in \partial X$ and  ${\bf y}$,  we will denote by ${\mathcal H}^{(a)}_x({\bf y}) $
the horoball (with respect to the metric $g_a$) centered at $x$ and passing through ${\bf y}$ ; furthermore,  for any ${\bf x} \in X$ we denote by $\psi_{x, {\bf y}}({\bf x}) $ its projection (with respect to $g_a$) on  the horosphere $\partial {\mathcal H}^{(a)}_x({\bf y})$.



{\bf In order to simplify the argument, one first assume that the two following conditions hold}
\begin{itemize}
\item (${\rm \bf C}_1$)   for any $x\in {\mathcal U}_p\cap \partial X$  the points  $h^n\cdot {\bf o}, n\in \Z^*, $ lie outside the horoball ${\mathcal H}_x^{(a)}(\bf o)$.
\item   (${\rm \bf C}_2$)   for any $x\in {\mathcal U}_h\cap \partial X$  the points  $p^m\cdot {\bf o}, m\in \Z^*, $ lie outside the horoball ${\mathcal H}_x^{(a)}(\bf o)$.
\end{itemize}

Fix  $k\geq 1$ and $\gamma \in \Gamma_{2k}$ with last letter in $h$. Let us decompose $\gamma$ into $a_{2k}a_{2k-1}\cdots a_1$ with 
$a_{2i}=p^{m_i}$ and $a_{2i-1}= h^{n_i}$ for $1\leq i\leq k$ ; set 
$\gamma_0:= Id$ and $\gamma_j:= a_j\cdots a_1$ for $1\leq j\leq 2k$.  We fix $x \in {\mathcal U}_p\cap \partial X$ ; by the ping-pong dynamic, there exists $c>0$ independent of $\gamma$ such that the distances   
$ d_a({\bf o},\psi_{x, {\bf o}}(\gamma^{-1}\cdot {\bf o}))
 $  and $ d_{a'}({\bf o},\psi_{x, {\bf o}}(\gamma^{-1}\cdot {\bf o}))
 $ are both $ \leq c$. 
 
 The cocycle property of the Busemann function thus  leads to the following
 \begin{eqnarray*} 
d_{a}({\bf o}, \gamma\cdot{\bf o}) 
&\geq& d_{a}({\bf o}, \psi_{x, {\bf o}}(\gamma^{-1}\cdot {\bf o}))-c\\
&=&   {\mathcal B}_{ x}^{(a)}(\gamma^{-1}\cdot {\bf o}, {\bf o})-c\\
 &=&  \sum_{j=0}^{2k-1}{\mathcal B}_{ x}^{(a)}(\gamma_{j+1}^{-1}\cdot {\bf o}, \gamma_j^{-1}\cdot {\bf o})-c\\
 &=&  \sum_{j=0}^{2k-1}{\mathcal B}_{\gamma_j\cdot x}^{(a)}(a_{j+1}^{-1}\cdot {\bf o}, {\bf o})-c,
  \end{eqnarray*}
and one may thus write, as in (\ref{defiruelle iterates})
\begin{equation}\label{ineq pour d_a}
\sum_{\stackrel {\gamma \in \Gamma_{2k}}{l(\gamma)=h}} e^{-{1\over 2}d_a({\bf x}, \gamma\cdot{\bf x})}  \leq e^{c\over 2} \mathcal L_a^{2k}1(x).
\end{equation}
By the previous assumption, all the quantities  ${\mathcal B}_{ x}^{(a)}(\gamma_{j+1}^{-1}\cdot {\bf o}, \gamma_j^{-1}\cdot {\bf o})$ above are non negative and 
  we want to compare them with a similar one involving $g_{a'}$. For any $x\in \partial X$ and ${\bf x, y} \in X$,  the quantity $\mathcal B_x({\bf x}, {\bf y})$ is equal to the ''signed'' length ( for $g_a$) of $[{\bf x}, \psi_{x, \bf y}({\bf x})]_a,$ the geodesic segment (for $g_a$)  joigning ${\bf x}$ and $\psi_{x, \bf y}(\bf x)$ ; in otherwords, with obvious notations, one gets 
  $$
  {\mathcal B}^{(a)}_x({\bf x}, {\bf y})= \int_{[{\bf x}, \psi_{x, \bf y}({\bf x})]_a}dg_a
  $$
  where the integral is non negative when ${\bf x}$ is outside   
${\mathcal H}_{ x}^{(a)} ({\bf y})$ and  negative when it lies inside. Similarly, we   introduce  the quantity $\beta_x({\bf x}, {\bf y})$ defined by
  $$
\beta_x({\bf x}, {\bf y})=  \beta_x^{(a, a')}({\bf x}, {\bf y}):= \int_{[{\bf x}, \psi_{x, \bf y}({\bf x})]_a}  dg_a'.$$
Note that for any ${\bf x, y,  z} $ in  $X$ and $\gamma \in \Gamma$ one gets  $  \beta_x({\bf x}, {\bf y})+  \beta_x({\bf y}, {\bf z})= \beta_x({\bf x}, {\bf z})$ and   $ \beta_x({\bf x}, {\bf y})=\beta_{\gamma \cdot x}(\gamma \cdot {\bf x}, \gamma \cdot {\bf y})$. 

 Since $ d_{a'}({\bf o},\psi_{x, {\bf o}}(\gamma^{-1}\cdot {\bf o}))
 $ is $\leq c$, we may write, as above  
\begin{eqnarray*} 
d_{a'}({\bf o}, \gamma\cdot{\bf o}) 
&\leq& d_{a'}({\bf o}, \psi_{x, {\bf o}}(\gamma^{-1}\cdot {\bf o}))+c\\
&\leq &
\beta_{x}(\gamma^{-1}\cdot {\bf o}, {\bf o})+c\\
 &=& \sum_{j=0}^{2k-1}\beta_{x}(\gamma_{j+1}^{-1}\cdot {\bf o}, \gamma_j^{-1}\cdot {\bf o})+c\\
 &=&  \sum_{j=0}^{2k-1}\beta_{\gamma_{j}\cdot x}(a_{j+1}^{-1}\cdot {\bf o}, {\bf o})+c.
  \end{eqnarray*}
which leads to the following inequality 
\begin{equation}\label{ineq pour d_a'}
\sum_{\stackrel {\gamma \in \Gamma_{2k}}{l(\gamma)=h}} e^{-{1\over 2}d_{a'}({\bf x}, \gamma\cdot{\bf x})}  \geq e^{-{c\over 2}} \mathcal K^{2k}1(x),
\end{equation}
 where $\displaystyle 
 {\mathcal K} \phi(y) := \sum_{\gamma \in \{p, h\}}\sum_{n \in 
 \Z^*}1_{x \notin U_{\gamma}}e^{-{1\over 2}\beta_{y}(\gamma^{-n}\cdot{\bf o} ,{\bf o})}\phi(\gamma^n\cdot y) $
  for any   function $\phi \in \mathbb L^{\infty}(\partial X)$ and any $y \in \partial X$.
 To prove (\ref{ineq-stricte}) it is thus sufficient to compare the spectral radius of $\mathcal L_a$ and $\mathcal K$ ; we will  use the following
  \begin{fact}\label{ineqdanshoroballe}
  For any $y \in \partial X$ and ${\bf x, y }\in X$   one gets 
   $$
  \Big\vert \beta_y({\bf x}, {\bf y}) \Big\vert  \leq 
    \Big\vert \mathcal B_y^{(a)}({\bf y}, {\bf y}) \Big\vert.
    $$
     Furthermore, for any $n \in \Z^*$, there exists $\eta(n)\geq 0$,  with $\eta(n)>0$ when $\vert n\vert $ is large enough, such that
    $$
\forall y\in {\mathcal U}_h \quad  
0\leq \beta_{y}(p^n\cdot {\bf o}, {\bf o})\leq  \mathcal B_y^{(a)}(p^n\cdot {\bf o}, {\bf o})-\eta(n).$$
  \end{fact}
\noindent Proof. The first large inequality is a direct consequence  of the   Remark  after Property  \ref{quasi-isometry}, namely $g_{a'}\leq g_a$. To prove the second one, we note that for any  $y \in \mathcal U_h$ and  any $ n\in \Z$ with $\vert n\vert$ large enough,  the geodesic segment 
$[p^{n}\cdot {\bf o}, \psi_{x, \bf o}(p^{n}\cdot {\bf o})]_a
$ inters sufficiently inside the horoball $\mathcal H$ centered at $\xi_p$ and in particular in the area where $g_a$ and $g_{a'}$ differ (ie $g_{a'}>g_a$) ; consequently $\beta_y(p^{n}\cdot {\bf o}, {\bf o}) -{\mathcal B}_y^{(a)}(p^{n}\cdot {\bf o}, {\bf o})>0$. the existence of $\eta(n)>0$ follows by an argument of continuity  with respect to $y$.\fdem

By this Fact, if $y \in \mathcal U_p$, one gets
$$
  \mathcal L_a1(y)=
 \sum_{n \in 
 \Z^*} e^{-{1\over 2}\mathcal B_{y}^{(a)}(h^{-n}\cdot{\bf o} ,{\bf o})}
 \leq \sum_{n \in 
 \Z^*} e^{-{1\over 2}\beta_{y}(h^{-n}\cdot{\bf o} ,{\bf o})} =\mathcal K1(y).
 $$
Assume now $y \in \mathcal U_h$ and fix $n_0\geq 1$ such that $\eta(n_0)>0$. By Property \ref{quasi-isometry},  one gets
$0\leq \beta_{y}(p^{-n_0}\cdot{\bf o} ,{\bf o})\leq K_0 d_0(p^{-n_0}\cdot{\bf o}, {\bf o})$ where $K_0$ is the constant which appears in $(\ref{K_0 quasi isometrie})$  ; consequently
$$
e^{-{1\over 2}\beta_{y}(p^{-n_0}\cdot{\bf o} ,{\bf o})}\geq \delta_0:=e^{-{K_0\over 2}d_0(p^{-n_0}\cdot{\bf o}, {\bf o})}.$$
On the other hand, by the above
  $$  \sum_{n \in 
 \Z^*} e^{-{1\over 2}\beta_{y}(p^{-n}\cdot{\bf o} ,{\bf o})}
\leq 
\sum_{n \in 
 \Z^*} e^{-{1\over 2}(d_{a'}(p^{-n}\cdot{\bf o} ,{\bf o})-c)}
\leq \Delta_0:=\sum_{n \in 
 \Z^*} e^{-{1\over 2K_0}(d_{0}(p^{-n}\cdot{\bf o} ,{\bf o})-c)}.
  $$  
 It follows
  \begin{eqnarray*}
   \mathcal L_a1(y)&=&e^{-{1\over 2}\mathcal B_{y}^{(a)}(p^{-n_0}\cdot{\bf o} ,{\bf o})}+ \sum_{\stackrel{n \in 
 \Z^*}{n\neq n_0}} e^{-{1\over 2}\mathcal B_{y}^{(a)}(p^{-n}\cdot{\bf o} ,{\bf o})}\\
 &\leq &e^{-{\eta(n_0)\over 2}}\times e^{-{1\over 2}\beta_{y}(p^{-n_0}\cdot{\bf o} ,{\bf o})}+ \sum_{\stackrel{n \in 
 \Z^*}{n\neq n_0}} e^{-{1\over 2}\beta_{y}(p^{-n}\cdot{\bf o} ,{\bf o})} \\
 &\leq& \rho \sum_{n \in 
 \Z^*} e^{-{1\over 2}\beta_{y}(p^{-n}\cdot{\bf o} ,{\bf o})} =\rho \mathcal K1(y),
  \end{eqnarray*}
  with $\rho:= 1-\Bigl(1-e^{-{\eta(n_0)\over 2}}\Bigr){\delta_0\over \Delta_0} \in ]0, 1[$.

Combining the two inequalities  $ \mathcal L_a1(y)\leq \mathcal K1(y)$ for $ y \in \mathcal U_p$ and $ \mathcal L_a1(y)\leq \rho \mathcal K1(y)$ for $ y \in \mathcal U_h$, one obtains by iteration
$$
\forall k \geq 1 \quad \mathcal L_a^{2k}1(.)\leq \rho^k \mathcal K^{2k}1(.)
$$
We put together this inequality with (\ref{ineq pour d_a}) and (\ref{ineq pour d_a'}) and obtain finally
$$
\sum_{\stackrel {\gamma \in \Gamma_{2k}}{l(\gamma)=h}} e^{-{1\over 2}d_a({\bf x}, \gamma\cdot{\bf x})}  \leq e^c \rho^k 
\sum_{\stackrel {\gamma \in \Gamma_{2k}}{l(\gamma)=h}} e^{-{1\over 2}d_a'({\bf x}, \gamma\cdot{\bf x})}.$$
  This gives the expected inequality $(\ref{ineq-stricte})$, in the case when 
   conditions $({\rm \bf C}_1)$ and $({\rm \bf C}_2)$ hold. 
   
When one or both of these  conditions do  not hold,  one replaces  the family $\{h^n: n \in \Z^*\} $ (resp.  $\{p^n: n \in \Z^*\} $) by the countable set 
 $H:= \{g \in \Gamma^{2N+1}/ l(g)= h\}$ (resp.    
 $P:= \{g \in \Gamma^{2N+1}/ l(g)= p\}$), where $N$ is choosen large enough such that   
  \begin{itemize}
\item    for any $x\in {\mathcal U}_p\cap \partial X$, the points  $ g \cdot {\bf o}, g \in H, $ lie outside the horoball ${\mathcal H}_x^{(a)}(\bf o)$.
\item for any $x\in {\mathcal U}_h\cap \partial X$, the points  $g \cdot {\bf o}, g \in P, $ lie outside the horoball ${\mathcal H}_x^{(a)}(\bf o)$.
\end{itemize}
Any $\gamma$ in $\Gamma_{2k(2N+1)}$  with last letter $h$ may be decomposed into 
 $\gamma = a_{2k}\cdots a_1$ with $a_{2i} \in P$ and $a_{2i-1} \in H$ for $1\leq i\leq k$ ; the same argument  as above, with obvious modifications, leads to the inequality
$$
\sum_{\stackrel {\gamma \in \Gamma_{2k(2N+1)}}{l(\gamma)=h}} e^{-{1\over 2}d_a({\bf x}, \gamma\cdot{\bf x})}  \leq e^c \rho^k 
\sum_{\stackrel {\gamma \in \Gamma_{2k(2N+1)}}{l(\gamma)=h}} e^{-{1\over 2}d_a'({\bf x}, \gamma\cdot{\bf x})}, $$
and   $(\ref{ineq-stricte})$ follows again.

\end{document}